\documentclass[24pt]{article}
\usepackage{amssymb}  
\usepackage{amsfonts}
\usepackage{amsmath}
\usepackage[T1]{fontenc}
\usepackage[italian]{babel}
\usepackage{pifont}
\usepackage{caption}
\usepackage{placeins}
\usepackage{pstricks}
\usepackage{pdfpages}
\usepackage{amsfonts}
\usepackage{graphicx}
\usepackage{placeins}
\usepackage{amsmath}
\usepackage{authblk}
\usepackage{pstricks}
\usepackage{pstricks-add}
\usepackage{framed}
\usepackage{bm}

\usepackage{tikz,fp,ifthen,fullpage}
\usepackage{pgfmath}
\usetikzlibrary{backgrounds}
\usetikzlibrary{decorations.pathmorphing,backgrounds,fit,calc,through}
\usetikzlibrary{arrows}
\usetikzlibrary{shapes,decorations,shadows}
\usetikzlibrary{fadings}
\usetikzlibrary{patterns}
\usetikzlibrary{mindmap}
\usetikzlibrary{decorations.text}
\usetikzlibrary{decorations.shapes}

\begin{document}

\bibliographystyle{plain}

\title{La costruzione di una scala musicale attraverso i numeri}

\author{F.~Talamucci}
\date{}

\maketitle

\section{Preambolo}

\noindent
Numerosi sono gli argomenti e i ragionamenti che, in modo pi\`u o meno spontaneo, mirano a fissare un'affinit\`a speciale tra la matematica e la musica: quest'ultima, pi\`u che altre attivit\`a artistiche,  sembra fortemente correlata con i numeri, con la geometria, con i calcolatori.

\noindent
Nel vasto repertorio di temi in proposito possiamo delineare due direzioni principali:

\begin{itemize}
\item[$(1)$] le regole della musica sono illustrate e decifrate dall'aritmetica, dalle leggi che regolano i numeri, 

\item[$(2)$] la struttura, la forma della musica viene resa comprensibile dalla geometria.
\end{itemize}

\noindent
Il secondo aspetto, senz'altro avvincente, pu\`o nascere dall'esigenza di ricerca di un sostegno, un appoggio razionale per dominare un'arte di per s\'e immateriale: la situazione non \`e molto diversa da quella di un pittore, uno scultore, un architetto che delineano la propria opera con le tracce 
esatte della geometria.

\noindent
In fondo, anche gli elementi pi\`u spontanei di un brano musicale, come il ritmo, una sequenza di accordi, 
lo sviluppo di una melodia, possono essere ricondotti in astratto a semplici figure geometriche.

\noindent
La perplessit\`a che talvolta emerge in relazione al punto $(2)$ \`e l'intenzionalit\`a, la coscienza da parte del compositore riguardo al disegno, allo schema che sembra affiorare dalla musica scritta: 
se per l'arte figurativa l'intenzione \`e quasi sempre molto chiara ed unanimamente raccolta, nel caso della musica spesso si forzano alcune componenti, alcuni spunti magari inconsapevoli; inoltre pu\`o apparire estremamente distante l'effettivo ruolo musicale di strutture astratte, spesso complesse, su cui un brano dovrebbe appoggiarsi.

\noindent
Il tema dapprima elencato come $(1)$ riunisce argomenti usuali e consolidati del binomio o matematica--musica: 
sin dai tempi antichi \`e apparso inevitabile l'impiego dell'artimetica nel percorso che stabilisce la gamma dei suoni utilizzabili e le relazioni tra essi. Si tratta di rintracciare la presenza della matematica non nella musica gi\`a formata, ma nel contesto delle regole musicali alla base del sistema adottato.

\noindent
Pi\`u che la conformazione si esplora dunque la formazione della musica, in un ambito che con maggiore precisione pu\`o essere attinente all'acustica, o alla psicoacustica\footnote{disciplina che studia la percezione soggettiva del suono, legata anche a fattori fisiologici e psicologici}, quando i numeri esigono di spiegare sensazioni come la consonanza.

\noindent
La nostra esposizione percorre nel contesto della tematica $(1)$ le tre modalit\`a,  probabilmente le pi\`u significative dal punto di vista storico, con le quali possiamo fissare i suoni di una scala musicale, che poi verranno combinati in orizzontale ed in verticale per produrre melodie e accordi.

\noindent
Il metodo proposto concede un privilegio al punto di vista della matematica, inseguendo l'obiettivo di 
for\-ma\-liz\-za\-re le tre scale musicali esaminate a partire da un'assioma, un \emph{principio}, dal quale discendono  corollari e propriet\`a: in tale sviluppo andr\`a riconosciuto il significato musicale dei vari passaggi logici e delle definizioni introdotte dal principio in avanti.
La risposta pi\`u semplice da parte della matematica, suoni equidistanti, d\`a origine alla scala pi\`u recente: questo \`e il motivo per cui dal punto di vista cronologico si opera un ri\-bal\-ta\-men\-to, cominciando dalla scala pi\`u recente, quella equabile; in effetti la matematica che permea tale sistema \`e la pi\`u spontanea, la pi\`u adatta a compiere operazioni sugli oggetti definiti. Tali operazioni corrispondono a procedure che comunemente compie anche chi ha un approccio dilettantesco con la musica, come ad esempio trasportare in acuto o in grave una melodia.

\noindent
I numeri della scala equabile sono gli irrazionali: questo sembra porre un ostacolo alla con\-si\-sten\-za musicale di tale scelta (e l'opposizione in tal senso fu un fatto storicamente avvenuto): vedremo come spiegazioni di tipo sperimentale vadano tuttavia a sorreggere la validit\`a del sistema musicale equabile, che \`e accolto praticamente dalla quasi totalit\`a della musica di oggi.

\noindent
Il secondo tema sviluppato \`e la costruzione della scala pitagorica partendo, anche qui, 
da un principio ed insistendo sugli aspetti matematici di scala infinita, di numeri razionali, ragionevolmente confrontati con i suoni della scala equabile.
Dal punto di vista del contenuto possiamo affermare che, se da una parte la costruzione dei suoni pitagorici 
percorre temi consueti dell'intervento matematico sulla musica, si \`e cercato dall'altra di rendere il pi\`u sistematico e logico possibile il metodo deduttivo, per ovviare alla difficolt\`a di alcune presentazioni che, di quando in quando, abbandonano il nesso matematico per asserzioni da accogliere come buone.

\noindent
Infine, la scala cosiddetta naturale viene costruita, anche in questo caso, a partire da una legge che stabilisca quali suoni \`e possibile includere. 
Rispetto allo sviluppo tradizionale dell'argomento, si \`e voluto insistere sulla possibilit\`a di codificare
i suoni della scala naturale poggiandosi esclusivamente sulla costruzione geometrica di media armonica, mostrando il carattere esclusivo della scelta operata dalla prassi musicale, se si adotta un criterio di selezione basato sulla media armonica. 
Questo approccio rappresenta, a nostro parere, una novit\`a rispetto al panorama dei metodi offerti.

\noindent
La proposta qui sviluppata si serve intenzionalmente un formalismo semplice che richiede conoscenze ma\-te\-ma\-ti\-che non certo impegnative: le nozioni di insieme dei numeri razionali e irrazionali, di proporzione sono sufficienti. Non c'\`e da aspettarsi una teoria elegante e sofisticata, bens\`i un percorso aritmetico rudimentale e tutto sommato abituale, che per\`o vuole contraddistingersi per l'autorit\`a della deduzione e della logica, fin tanto che sia possibile, sull'acquisito e il non spiegato, e per l'intervento di locuzioni della teoria musicale e di rimandi di tipo storico solo se ritenuto inevitabile.
Dal punto di vista della musica, si \`e cercato di rendere il meno possibile indispensabile la conoscenza 
delle primarie nozioni musicali (scala, note, intervallo, ...), nell'intenzione (senz'altro ambiziosa) di accostarsi a tali concetti attraverso la matematica al tempo stesso contando sull'indubitabile esperienza di aver osservato qualche strumento musicale, di aver ascoltato musica in genere, ..
Per quanto riguarda l'acustica, \`e sufficiente la nozione di frequenza di un suono (come numero di vibrazioni in un secondo che danno origine al suono) e della propriet\`a ad essa collegata di suono pi\`u acuto o suono pi\`u grave, misurato dall'aumento o dalla diminuzione, rispettivamente, della frequenza. Inoltre, la conoscenza della legge (cosiddetta pitagorica) di proporzionalit\`a inversa tra la frequenza del suono e la lunghezza della corda in vibrazione che lo produce \`e senz'altro auspicabile.
La nozione, infine, di suono complesso formato dalla sovrapposizione di suoni semplici con frequenze doppia, tripla, quadrupla, ... pu\`o aiutare a comprendere alcune osservazioni, comunque a carattere secondario.

\section{La scala musicale}

\noindent
La musica adopera i suoni per formare una {\it melodia}, ovvero una successione di suoni, oppure un {\it accordo}, ovvero una sovrapposizione di due o pi\`u suoni. L'{\it armonia musicale} \`e l'insieme delle regole atte a formare e studiare gli accordi.

\noindent
Poniamoci la domanda: quali suoni posso utilizzare, ovvero quali sono le frequenze dei suoni che entrano a far parte dell'alfabeto musicale? 

\noindent
Se pensiamo alla forma pi\`u naturale di produrre musica, attraverso il {\it canto}, non sembrano esserci esigenze di questo tipo: nell'ambito dell'estensione della propria voce tra le frequenze $f_m$ e $f_M$ posso emettere il suono a qualunque frequenza $f$ compresa tra $f_m$ $f_M$.  

\noindent
Alcuni strumenti offrono questa medesima possibilit\`a: sul violino, strisciando lungo la tastiera, il suono ottenuto varia con continuit\`a tra il suono iniziale e quello finale ({\it glissando}). Anche il trombone \`e uno strumento che permette la continuit\`a tra un suono e l'altro.

\noindent
Tuttavia, la maggior parte degli strumenti musicali (tutti gli strumenti a tastiera, strumenti a corda come la chitarra, l'arpa, ..., strumenti a fiato con i fori o con le {\it chiavi}, ovvero i tasti da premere per chiudere i fori, strumenti a percussione come il timpano, lo xilofono, ...), prevedono solo la possibilit\`a di produrre un numero finito di frequenze: i tasti del pianoforte, premuti dal pi\`u grave al pi\`u acuto, danno luogo a $88$ suoni, ovvero $88$ frequenze diverse, in ordine crescente.
La corda di una chitarra fornisce solamente i suoni che ottengo premendo gli spazi sulla tastiera separati dalle sbarrette trasversali in rilievo ({\it tasti}). La frequenza del timpano viene modificata, attraverso un pedale, portando la tensione della membrana ad un valore maggiore o minore: in generale, solo poche frequenze sono producibili.

\noindent
Tutti questi strumenti in grado di emettere solo particolari frequenze vengono detti {\it a suono determinato}. 

\noindent
A questo punto \`e chiaro che bisogna stabilire l'altezza dei suoni ai quali tutti gli strumenti (a suono determinato  o no) devono riferirsi ed adeguarsi per eseguire musica. 

\noindent
Diciamo che la gamma dei suoni fissati, indipendentemente dal particolare strumento, sia formata dalla sequenza crescente degli $N$ suoni $\{ f1, f2, \dots, fN\}$: adopereremo per essa il termine {\it scala}.
''Accordare'' significa scegliere un particolare insieme di frequenze (per le corde di uno strumento, per l'intonazione di uno strumento a fiato, ...). Il problema dell'accordatura, in questo senso, diventa quello di fissare una scala.

\subsection{Il Principio dell'ottava}

\noindent
Il primo aspetto che vogliamo analizzare consiste nello stabilire il numero minimo di frequenze da fissare per avere il panorama completo dei suoni consentiti: ovvero, vanno elencati tutti, nell'arco dell'intera gamma delle frequenze udibili (all'incirca dai $16$, $20$ $Hz$ fino a $20.000$ $Hz$) oppure qualche argomento di natura musicale comporta una sorta di ripetizione della regola di formazione della scala, replica operata su determinati segmenti della gamma sonora?

\noindent
Il primo argomento che fa da cardine alla costruzione dei vari tipi di scale musicali viene qui indicato 
come

\noindent
{\bf Principio dell'Ottava}: se una scala contiene il suono $f$, contiene anche i suoni 
$2f$, $4f$, $8f$, $\dots$, $2^Nf$, $\dots$ e i suoni 
$\frac{1}{2}f$, $\frac{1}{4} f$, $\frac{1}{8}f$, $\dots$, $\frac{1}{2^N}$, $\dots$
per ogni $N=1,2,3,\dots$. 
Sinteticamente, possiamo affermare che la presenza di $f$ nella scala fa includere nella medesima scala tutti i suoni 
\begin{equation}
\label{2kf}
2^kf, \quad k\in {\Bbb Z} \quad \textrm{numeri interi}.
\end{equation} 

\noindent
Un Principio \`e un postulato che si accetta, si ammette valido {\it a priori} per portare avanti una teoria, sviluppando le conseguenze; sarebbe tuttavia scorretto non evidenziare le ragioni di natura musicale che giustificano il Principio. 

\noindent
Si pu\`o effettuare  un esperimento, davanti ad una platea di individui dotati di differenti estensioni di voce. Alla richiesta  di intonare un suono abbastanza acuto o abbastanza grave, \`e naturale riscontrare che, nel primo caso, chi non arriva a tale altezza (tipo una voce maschile) produrr\`a istintivamente il suono proprio ad una frequenza dimezzata, o divisa per quattro, ... . ; 
nel secondo caso, una voce femminile non provvista di suoni gravi intoner\`a il suono a frequenza doppia, o quadrupla, ... . In una situazione simile a questa, 
ciascun corista di un coro  ``capta'' i $440$ $Hz$ del diapason, messo in vibrazione prima dell'esecuzione per intonarsi, assimilandolo secondo la propria estensione. 

\noindent
Allo stesso modo, prima di un'esecuzione orchestrale ogni strumento di un'or\-che\-stra si adegua al LA 
dell'accordatura, ciascuno secondo la propria estensione.
Dunque, il suono non riproducibile dalla propria gamma sonora
viene istintivamente giudicato ``sostituibile'' dai suoni $2^kN$, $k\in {\Bbb Z}$ indicati dal Principio. 

\noindent
La naturale spontaneit\`a del Principio \`e avvalorata dalla consuetudine di \emph{assegnare il medesimo nome per i suoni equivalenti}: l'identica denominazione da parte della teo\-ria musicale delle note con frequenze come in (\ref{2kf}) sottolinea la forte somiglianza ed omogeneit\`a fra tali suoni. Per fare un esempio, oltre al LA centrale a $440$ $Hz$, sul pianoforte sono presenti altri sette LA, alle frequenze (in $Hz$)
$27,5$, $55$, $110$, $220$, $880$, $1760$, $3520$.

\noindent
Il nome \emph{ottava} proviene dal fatto che, secondo la successione dei sette suoni DO, RE, MI, FA, SOL, LA SI, DO formanti la cosiddetta \emph{scala diatonica} si ritrova il medesimo dopo sette, cio\`e all'ottava nota. Il DO alla fine di tale sequenza ha frequenza doppia rispetto al DO iniziale. 

\noindent
L'accoglimento del Principio dell'ottava \`e un fatto indiscusso in ogni civilt\`a musicale. 
Gi\`a nell'Antica Grecia, l'ottava veniva indicata con il termine $\delta \iota \alpha$ $\pi \alpha \sigma \omega \nu$ , diapason, ovvero che passa attraverso ({\it di\`a}) tutti ({\it pas\`on}) i toni, prima di ritrovarne uno ``identico'' a quello di partenza.

\noindent
Oltre alla concreta evidenza della validit\`a del Principio attraverso i fatti che abbiamo descritto, 
si possono addurre ragioni di natura differente:

\begin{itemize}
\item[$(i)$] produrre un suono a frequenza doppia richiede l'operazione pi\`u semplice che si possa effettuare sull'oggetto che emette il suono: la corda divisa a met\`a, la canna divisa a met\`a, ... La spontaneit\`a con cui un mezzo musicale produce due suoni identici mediante la divisione per due, d\`a luogo ad una motivazione legata, per cos\`\i$\,$dire, al numero, al rapporto semplice: l'origine del presupposto va  rintracciata nelle idee pitagoriche circa la Musica regolata, come tutto l'Universo, dal Numero. In quest'ordine di idee, 
nel Principio matematica e musica partecipano equamente. Per secoli, \`e prevalsa l'idea
della Musica come scienza esatta che si appoggia all'aritmetica: in questo senso, la semplice divisione per due deve dare origine, in Musica, ad un effetto inevitabimente favorevole.

\item[$(ii)$] La prima tra le frequenze delle armoniche di una frequenza fondamentale \`e proprio quella a frequenza doppia: la modalit\`a di vibrazione della corda suddividendosi in due parti uguali \`e la prima che si incontra tra quelle consentite. In altre parole, il primo suono tra quelli che vanno a formare un suono composto di frequenza fondamentale $f$ \`e quello a frequenza $2f$.

\item[$(iii)$] Si pu\`o aggiungere una motivazione di tipo fisiologico, ricordando che suoni di differente frequenza mettono in moto nell'orecchio regioni diverse dell'organo di ricezione, che ha la forma di un'elica. Le lunghezze delle porzioni di elica che competono ad una specifica ottava sono pressoch\'e costanti, ovvero i punti dell'elica in cui risuonano e sono percepite le frequenze (\ref{2kf}) (nell'ambito del range di udibilit\`a) risultano equidistanti.
\end{itemize}

\noindent
Il Principio dell'ottava permette di semplificare la definizione di Scala musicale: in effetti, 
richiedere la presenza di un suono $f$ in tutte le ottave che mi occorrono fissa una modalit\`a operativa per le frequenze da fissare nella gamma dei suoni.
\noindent
Pi\`u precisamente, lo stabilire gli $N-1$ suoni intermedi $f_1$, $\dots$, $f_{N-1}$ nel dominio dell'ottava $\left( f_0, 2f_0\right)$ d\`a luogo alla scala delle $N+1$ frequenze crescenti
\begin{equation}
\label{scala}
f_0, \;f_1, \;f_2,\;\dots\;,\; f_{N-1}, \;f_N=2f_0.
\end{equation}
La (\ref{scala}) comporta automaticamente, in base a (\ref{2kf}), la suddivisione delle altre ottave secondo
\begin{equation}
\label{altrescale}
2^k f_0, \quad 2^k f_1, \quad 2^k  f_2, \;\dots\;,\quad\;2^k f_N=2^{k+1}f_0 
\end{equation}
dove $k=0$ denota la scala nell'ottava di riferimento, $k=1,2,\dots$  la scala nelle ottave pi\`u acute,
$k=-1, -2, \dots$ in quelle pi\`u gravi, ovviamente nell'ambito dell'estensione dello strumento o, pi\`u in generale, dell'udibilit\`a. 
Possiamo pertanto perfezionare la definizione di \emph{scala musicale} affermando che essa \`e la  suddivisione dell'ottava in un numero prestabilito di suoni di frequenza intermedia, che sono a disposizione per essere combinati in sequenza (melodie) o in sovrapposizione (accordi). La gamma dei suoni intermedi in un'ottava svolge il ruolo di ``alfabeto'' dei suoni per comporre ed eseguire musica.
Indichiamo la scala compresa tra i suoni $f_0$ e $2f_0$ mediante $[f_0, 2f_0]$. 

\noindent
Andiamo ora a compiere un passo logico importante: il principio assoluto accettato da ogni civilt\`a musicale, cui ci siamo riferiti come principio dell'ottava, comporta, per cos\`\i$\,$dire, una \emph{metrica}, una misura della distanza tra suoni.

\subsection{La ``distanza'' sulla scala: gli intervalli}

\noindent
Pensiamo alla totalit\`a della gamma delle frequenze come ad una linea continua e fissiamo su di essa $f_0$, frequenza di riferimento: in base a questa abbiamo individuato i suoni $2^kf_0$, con $k$ intero, alle ottave sopra e sotto: 
{\footnotesize 
$$
\begin{array}{l}
\; |-|--|----|--------|---------------|\\
\frac{f_0}{4} \;\, \frac{f_0}{2}\;\;\;f_0\quad \quad  \quad \;\;\;2f_0 \qquad \qquad \qquad 
\qquad \qquad \qquad \qquad \qquad \qquad \qquad\;
\quad \quad 8f_0
\end{array}
$$
}
\`E evidente che gli intervalli di frequenze che stiamo considerando uguali dal punto di vista della costruzione della scala hanno lunghezze ``euclidee'' differenti: la lunghezza del segmento tra $2^kf_0$
e $2^{k+1}f_0$ che copre una singola ottava \`e  $2^{k+1}f_0 - 2^kf_0=2^kf_0$.
D'altra parte, l'identificazione dei suoni $2^kf_0$ deve generare una regola che permetta di ritenere ``uguali'' i segmenti con gli estremi indicati: per questo basta semplicemente osservare che 
il \emph{rapporto} tra i valori delle frequenze agli estremi rimane costante, indipendentemente da $k$:
\begin{equation}
\label{estremi}
\dfrac{2^{k+1} f_0}{2^kf_0}=2\quad \textrm{per ogni}\;k\in {\Bbb Z}
\end{equation}

\noindent
L'effetto  di rendere ``uguali'' i segmenti che contengono un'ottava valutando il rapporto anzich\'e la lunghezza effettiva rimane valido anche per gli altri suoni intermedi della scala (\ref{scala}), da confrontare con quelli presenti sulle scale indotte (\ref{altrescale}): 
al suono $f_{i_1}$ della scala $[f_0, 2f_0]$ corrisponde $2^kf_{i_1}$, per un $k$ fissato, della scala $[2^kf_0, 2^{k+1} f_0]$, medesima situazione per $f_{i_2}$ e $2^kf_{i_2}$. 
Avendo sulle rispettive scale il medesimo ruolo, le due coppie $(f_{i_1}, f_{i_2})$ e 
$(2^kf_{i_1}, 2^kf_{i_2})$ devono rappresentare segmenti di uguale lunghezza e ci\`o si realizza 
misurando il rapporto, dato che evidentemente
$\frac{2^kf_{i_1}}{2^kf_{i_2}}=\frac{f_{i_1}}{f_{i_2}}$ qualunque sia $k\in {\Bbb Z}$. 
La validit\`a del Principio dell'ottava comporta l'uguaglianza dei segmenti del tipo indicato nello schema qui di seguito:
$$
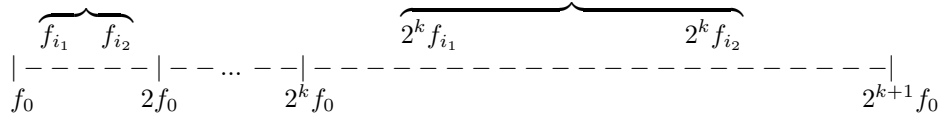
\begin{array}{l}
\;\;\;\;\overbrace{f_{i_1}\;\;\;\;f_{i_2}} \hspace{3.5truecm}\overbrace{2^kf_{i_1}\hspace{3truecm} 2^kf_{i_2}} \\
|-----|--\,...\,--|----------------------|\\
f_0 \qquad \qquad 2f_0 \qquad \qquad 2^kf_0 \hspace{7truecm} 2^{k+1}f_0
\end{array}
$$
La conclusione da trarre \`e la seguente:  
nella scala delle frequenze la grandezza significativa per misurare le distanze \`e il \emph{rapporto}: 
attraverso questa associazione, i suoni che si corrispondono su pi\`u scale hanno medesima distanza.
Definiamo \emph{distanza tra due suoni} $f_{I_1}$ e $f_{i_2}$ il rapporto 
\begin{equation}
\label{interv}
I(f_{i_1}, f_{i_2})=\dfrac{f_{i_2}}{f_{i_1}}, \qquad f_{i_2}\geq f_{i_1}.
\end{equation}
Tale definizione accoglie la richiesta $I(2^k f_{i_1}, 2^kf_{i_2})=I(f_{i_1}, f_{i_2})$ per ogni $k\in {\Bbb Z}$.

\noindent La lettera $I$ non \`e utilizzata a caso: infatti, in termini musicali la distanza (\ref{interv}) si chiama \emph{intervallo musicale}. Ciascun intervallo acquista un nome specifico; per il momento conosciamo i termini \emph{unisono}: $I(f, f)=1$ e \emph{ottava}: $I(f, 2f)=2$, qualunque sia il suono $f$. 

\noindent
Nell'ambito di una scala $[f_0, 2f_0]$, il valore di $I$ \`e sempre compreso tra l'unisono e l'ottava, nel sendo che comunque si scelgano $2$ elementi $f_{i_1}$, $f_{i_2}\geq f_{i_1}$ della scala $[f_0, 2f_0]$ 
si ha $1\leq I(f_{i_1}, f_{i_2}) \leq 2$. In tal modo, la scala $[f_0, 2f_0]$ \`e completamente individuata dagli $N$ numeri compresi tra $1$ e $2$
\begin{equation}
\label{gradi}
G_1=I(f_0, f_1)=\dfrac{f_1}{f_0}, \;\;\; G_2=I(f_1, f_2)= \dfrac{f_2}{f_1}, \; \;
\dots\;\; G_N=I(f_{N-1}, 2f_0)=\dfrac{2f_0}{f_{N-1}}
\end{equation}
che rappresentano le distanze tra gli elementi consecutivi.
Fissato $f_0$, dunque anche $f_N=2f_0$, i numeri effettivamente da conoscere per ricavare $G_1$, $\dots$, $G_N$ sono $N-1$: $f_1$, $\dots$, $f_{N-1}$.
Chiamiamo $G_1$, $G_2$, $\dots$, $G_N$ i {\it gradini della scala}.
La scala $[2^kf_0, 2^{k+1}f_0]$ \`e divisa mediante gradini identici, dato che 
$I(f_0, f_1)=I(2^k f_0, 2^k f_1)$, $I(f_1, f_2)=I(2^k f_1, 2^k f_2)$, $\dots$, 
$I(f_{N-1}, f_N)=I(2^k f_{N-1}, 2^k f_N)$
qualunque sia $k\in {\Bbb Z}$. 

\noindent
L'indagine aritmetica che ha condotto ad una misura relativa tra suoni necessita a questo punto di 
qualche considerazione e valutazione di tipo musicale. In effetti, 
dal punto di vista della percezione, la capacit\`a di riconoscimento di un brano musicale \`e proprio basata sul rapporto, ovvero sull'intervallo. 
Racchiude in s\'e maggior valore musicale il rapporto (intervallo) del suono assoluto (intonazione).
La riconoscibilit\`a avviene attraverso gli intervalli. 
Anche dal punto di vista storico, ogni sede aveva un proprio diapason di riferimento, pur mostrando i medesimi sviluppi musicali, le medesime modalit\`a nel combinare i suoni, nel percorrere i rapporti tra suoni.

\noindent
Una semplice esperienza in tal senso pu\`o essere quella di sperimentare quanto segue.
Con una chitarra  si pu\`o eseguire una successione di suoni sulla VI corda (suono $f_0$), quella pi\`u grave, posta pi\`u in alto di tutte, premendo il dito su una precisa sequenza di spazi. 
Trasferendo la medesima sequenza di spazi sulla prima corda (quella pi\`u in basso, con suono $4f_0$) si ascolta la medesima melodia, anche se pi\`u acuta. La medesima posizione in verticale degli spazi mantiene infatti inalterato il rapporto dei suoni corrispondenti.

\noindent
Una prova simile pu\`o essere effettuata su una tastiera: a partire da un qualunque tasto $f_0$ si riconosce la medesima sequenza di tasti bianchi o tasti neri dopo il dodicesimo tasto: ad esempio pu\`o essere 

$$
\begin{array}{lllllllllllll}
\square &\blacksquare  &\square & 
\square & \blacksquare  &\square  &
\blacksquare & \square &\blacksquare & 
\square & \square  &\blacksquare & \square  \\
f_0& f_1&f_2&f_3 &f_4&f_5&f_6&f_7&f_8&f_9&f_{10}&f_{11}& f_{12}=2f_0 \\
\\
\\
\square &\blacksquare  &\square & 
\square & \blacksquare  &\square  &
\blacksquare & \square &\blacksquare & 
\square & \square  &\blacksquare & \square  \\
2f_0& 2f_1&2f_2&2f_3 &2f_4&2f_5&2f_6&2f_7&2f_8&2f_9&2f_{10}&2f_{11}& 2f_{12}=4f_0 
\end{array}
$$
Premendo una qualsiasi fila di tasti nell'ottava $[f_0, 2f_0]$ e la sequenza nel medesimo ordine nell'ottava $[2f_0, 4f_0]$, si riconosce senz'altro la medesima melodia, una tra le quali pu\`o essere
$$
f_4 \;f_4\; f_5\; f_7\; f_7\; f_5\; f_4\; f_2\; f_0\; f_0\; f_2\; f_4\; f_4\;\;\;\; f_2.
$$

\subsection{Operazioni con gli intervalli}

\noindent
Fare riferimento al rapporto tra gli estremi induce la seguente propriet\`a di composizione degli intervalli:

\noindent
Se $I(f_{i_1}, f_{i_2})$ \`e un intervallo che comprende il suono intermedio $f_{i_3}$, ovvero 
$f_{i_1}\leq f_{i_2}\leq f_{i_3}$,  
allora 
\begin{equation}
\label{sommainterv}
I(f_{i_1}, f_{i_2})=I(f_{i_1}, f_{i_3})I(f_{i_3}, f_{i_2}).
\end{equation}
Infatti $I(f_{i_1}, f_{i_3})I(f_{i_3}, f_{i_2})=\dfrac{f_{i_3}}{f_{i_1}}\dfrac{f_{i_2}}{f_{i_3}}=
\dfrac{f_{i_2}}{f_{i_1}}=I(f_{i_1}, f_{i_2})$. 

\noindent
Rappresentando graficamente gli intervalli come segmenti, si evidenzia come il segmento 
somma di due segmenti adiacenti (in senso euclideo) ha come lunghezza il prodotto delle lunghezze dei due contributi:
$$
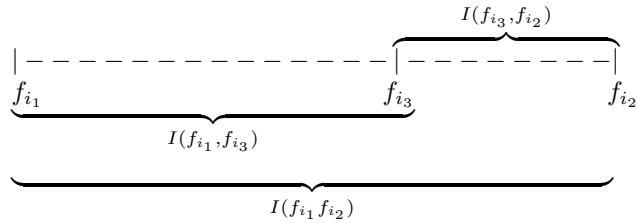
\begin{array}{l}
\quad |--------------\overbrace{|--------|}^{I(f_{i_3}, f_{i_2})}\\
\quad \underbrace{
f_{i_1} \qquad  \qquad \qquad \qquad \qquad \qquad \;\;\;f_{i_3}}_{I(f_{i_1}, f_{i_3})}\qquad \qquad \qquad \;\;\;\;f_{i_2}\\
\quad \underbrace{\hspace{8truecm}}_{I(f_{i_1} f_{i_2})}
\end{array}
$$
Questa operazione pu\`o anche uscire dall'ambito dell'ottava: ad esempio unire $2$ ottave con $3$ ottave d\`a
$I(f_0, 4f_0)I(4f_0, 32f_0)=I(f_0, 32f_0)=I(f_0, 2^5 f_0)$, ovvero $5$ ottave.

\noindent
Elaboriamo ora un concetto di uguaglianza tra insiemi di suoni.
Due figure geometriche sono congruenti se perfettamente sovrapponibili ovvero se coincidono a meno di isometrie. Esportando questa nozione a due insiemi di suoni, un'isometria, in quanto  
trasformazione da un insieme all'altro che non altera le distanze reciproche fra i suoni, dovr\`a mantenere inalterati i rapporti fra i suoni corrispondenti, avendo in mente che \`e la (\ref{interv}) ad esprimere la distanza fra due suoni. \`E in questo modo che formuliamo il concetto di uguaglianza tra due insiemi di suoni: pi\`u precisamente, l'insieme ordinato di suoni (non necessariamente appartenenti alla medesima ottava)
\begin{equation}
\label{seq}
\left\{ f_{i_1}, \;f_{i_2},\; \dots,\; f_{i_n}\right\}
\end{equation}
\`e \emph{congruente} all'insieme ordinato $\left\{ {\widetilde f}_{i_1}, \;{\widetilde f}_{i_2},\; \dots,\; {\widetilde f}_{i_n}\right\}$
contenente il me\-de\-si\-mo numero di suoni se i rapporti tra i suoni corrispondenti, nel medesimo ordine, sono i medesimi:

\begin{equation}
\label{congruenti}
\dfrac{{\widetilde f}_{i_1}}{f_{i_1}}=\dfrac{{\widetilde f}_{i_2}}{f_{i_2}}=
\dfrac{{\widetilde f}_{i_3}}{f_{i_3}}=
\;\dots \;=\dfrac{{\widetilde f}_{i_{n-1}}}{ f_{i_{n-1}}}
=\dfrac{{\widetilde f}_{i_n}}{ f_{i_n}}
\end{equation}
In modo equivalente, possiamo definire congruenti due insiemi ordinati di suoni che si sviluppano secondo i medesimi rapporti tra frequenze consecutive del medesimo insieme:
\begin{equation}
\label{congruenti2}
\begin{array}{lllll}
\dfrac{f_{i_2}}{f_{i_1}}=\dfrac{{\widetilde f}_{i_2}}{{\widetilde f}_{i_1}}, &
\dfrac{f_{i_3}}{f_{i_2}}=\dfrac{{\widetilde f}_{i_3}}{{\widetilde f}_{i_2}}, &
\dots &
\dfrac{f_{i_n}}{f_{i_{n-1}}}=\dfrac{{\widetilde f}_{i_n}}{{\widetilde f}_{i_{n-1}}}.
\end{array}
\end{equation}
Ad esempio, una sequenza che si articola sulle frequenze multiple ${\widetilde f}_{i_1}=\kappa f_{i_1}$, $\dots$, ${\widetilde f}_{i_n}=\kappa f_{i_n}$ \`e congruente a quella di partenza per qualunque $\kappa \in {\Bbb Z}$. 
Utilizzeremo anche la scrittura $\left\{ {\widetilde f}_{i_1},\dots,\; {\widetilde f}_{i_n}\right\}\equiv  \left\{ f_{i_1}, \; \dots,\; f_{i_n}\right\}$ per due insiemi di suoni congruenti.

\noindent 
Adoperando la definizione di congruenza (\ref{congruenti2}) per insiemi (\ref{seq}) che sono essi stessi scale musicali, troviamo che due scale $\{ f_1, \dots, f_N\}$ e $\{ {\widetilde f}_1, \dots, {\widetilde f}_N\}$ sono congruenti se si sviluppano esattamente sui medesimi gradini (\ref{gradi}):

$$
\begin{array}{lllll}
\dfrac{f_2}{f_1}=\dfrac{{\widetilde f}_2}{{\widetilde f}_1}, &
\dfrac{f_3}{f_2}=\dfrac{{\widetilde f}_3}{{\widetilde f}_2}, &
\dots &
\dfrac{f_n}{f_{n-1}}=\dfrac{{\widetilde f}_n}{{\widetilde f}_{n-1}}.
\end{array}
$$

\subsection{Vari tipi di scala}

\noindent
Sulla base comune del Principio dell'ottava si diramano, a seconda della scelta di $N$ e dei gradini (\ref{gradi}) le varie possibilit\`a di formare una scala. Ci si pu\`o basare, anche se in modo grossolano, sulle seguenti indicazioni:
\begin{itemize}
\item[$(A)$] $N$ \`e stabilito da una cultura musicale, o da un genere musicale: da secoli la nostra civilt\`a utilizza $N=12$, ovvero dopo $12$ passi da $f_0$ si trova $2f_0$, con $11$ suoni intermedi. Questo si pu\`o riscontrare sul pianoforte, dove a partire da un qualunque tasto, dopo dodici tasti (in s\`u o in gi\`u, tasti neri compresi) trovo il suono di partenza all'ottava sopra, o all'ottava sotto, oppure sulla chitarra: a partire da uno spazio, contando gli spazi lungo la medesima corda verso il ponte, si trova al dodicesimo spazio l'ottava.
La scala con $12$ suoni viene detta \emph{scala cromatica}; la \emph{scala dia\-to\-ni\-ca}, composta dai sette suoni tradizionalmente noti (da DO a SI) \`e un sottoinsieme di quella cromatica.

\item[$(B)$] Nell'ambito della scelta di $N$, fissare i gradini $G_1$, $\dots$, $G_N$ significa {\it accordare} la scala. Per inquadrare la questione nella nostra musica occidentale attraverso i secoli, il problema ha una dimensione, diciamo cos\`\i$\,$, storica: in alcune fasi \`e prevalsa la convinzione che le frequenze intermedie debbano provenire solo da calcoli ``semplici'' sul suono $f_0$, ovvero attraverso moltiplicazioni e divisioni con numeri interi. 

\noindent
In altre fasi, questo sistema \`e andato a scontrarsi con la difficolt\`a di poter essere adoperato, ad esempio, dagli strumenti a tastiera. 
La soluzione e\-sco\-gi\-ta\-ta propone una lista $G_1$, $\dots$, $G_N$ molto prossima a quella abbandonata che permette il superamento di alcune complicazioni.
\end{itemize}

\noindent
Non ci sono argomenti spontanei di tipo matematico che stabiliscano quanti suoni scegliere (fase $(A)$): la fase $(A)$ \`e propria della musica. In effetti, $N=12$ non ha una valenza universale e 
altri esempi possono essere raccolti da una particolare civilt\`a musicale, o della caratterizzazione di un genere musicale, o dalla ideazione da parte di un compositore: 
\begin{description}
\item[$N=5$]: scala pentatonica, utilizzata in Cina, anche dal Jazz,
\item[$N=6$]: scala esatonale, utilizzata in Oriente e, per esempio, da Debussy, 
\item[$N=7$]: scala blues, scala araba, scala di Bach, scala di Skrjabin, scala enigmatica,...: queste scale possono essere, sostanzialmente, comprese apportando qualche modifica alla scala diatonica,
\item[$N=8$]: scala di genere \emph{bebop}.
\end{description}
Dall'altra parte, si intuisce che il punto $(B)$, ovvero decidere {\it quali} suoni, ha un collegamento naturale con la matematica: storicamente, l'impostazione aritmetica ha portato a considerare la questione musicale nel campo della scienza dei numeri. 
In modo sommario, possiamo anticipare che in termini puramente aritmetici affrontare il problema $(B)$ significa porsi sostanzialmente di fronte alla scelta \emph{numeri razionali} oppure \emph{numeri irrazionali}.

\noindent
Le tre Sezioni che seguono presentano tre differenti scale musicali, le pi\`u importanti dal punto di vista storico, analizzate da un punto di vista aritmetico, ponendosi innanzitutto la questione di quali assiomi, oltre al principio dell'ottava, rendono definita la scala e rilevando poi le specifiche propriet\`a di ciascuna costruzione di suoni. 

\section{La scala equabile}

\noindent
Dal punto di vista della matematica, ignorando l'evoluzione storica del travagliato e controverso problema 
di fissare la scala dei suoni, partiremo dalla soluzione pi\`u semplice, quella per cui si vuole rendere la scala pi\`u omogenea possibile, nel senso che andremo a specificare. L'aggettivo stesso che qualifica tale scala, {\it equabile}, rende l'idea delle particolarit\`a che andiamo cercando.

\noindent
Poniamo attenzione alla (\ref{interv}): 
una scala equidistanziata, che chiameremo scala \fbox{E}, vuole che i gradini $G_1$, $\dots$, $G_N$ di (\ref{gradi}) siano tutti uguali e pari ad un numero positivo $r$: 
$$
G_1=G_2=\dots=G_N=r
\quad
\textrm{ovvero}
\quad
\dfrac{f_1}{f_0}=\dfrac{f_2}{f_1}=\dots = \dfrac{2f_0}{f_{N-1}}=r.
$$
Questo comporta che i suoni della scala formino una \emph{progressione geometrica} di ragione $r$. 
Per determinarli, operiamo ad esempio come segue:
$$
2=\dfrac{2f_0}{f_0}=\dfrac{2f_0}{f_{N-1}} \dfrac{f_{N-1}}{f_{N-2}}\dots \dfrac{f_2}{f_1}\dfrac{f_1}{f_0}=r^N
$$
da cui $r=\sqrt[N]{2}$. Per ottenere i valori dei suoni della scala \fbox{E} possiamo operare per mezzo della regola (\ref{sommainterv}) a partire da $\dfrac{f_1}{f_0}=r=\sqrt[N]{2}$:
$$
\dfrac{f_2}{f_0}=\dfrac{f_2}{f_1}\dfrac{f_1}{f_0}=r\cdot r=\left(\sqrt[N]{2}\right)^2=\sqrt[N]{2^2}, \quad
\dfrac{f_3}{f_0}=r^3=\sqrt[N]{2^3}, \quad \dots,\quad \dfrac{f_{N_1}}{f_{N-2}}=r^{N-1}=\sqrt[N]{2^{N-1}}
$$
e infine, come conferma, $\dfrac{2f_0}{f_0}=r^N=\sqrt[N]{2^N}=2$. Pertanto i suoni della scala \fbox{E} che parte da $f_0$ sono
\begin{equation}
\label{ognifreq}
f_K=(\sqrt[N]{2})^Kf_0, \qquad  K=0,1,\dots, N
\end{equation}
e per gli altri valori di $K\in {\Bbb Z}$ si ottengono i suoni disposti sulle altre ottave.

\noindent
Per $N=2$ e $N=3$, casi per\`o privi di interesse dal punto di vista musicale, si ottengono rispettivamente le scale $\{f_0, \sqrt{2}f_0, 2f_0 \}$, $\{f_0, \sqrt[3]{2}f_0,\sqrt[3]{2^2}, 2f_0\}$. 
I numeri $\sqrt{2}$, $\sqrt[3]{2}$, $\sqrt[3]{2^2}$ che formano le frequenze intermedie sono numeri irrazionali. La circostanza rimane vera in generale per $N>3$: i numeri $\sqrt[N]{2}$, $\sqrt[N]{2^2}$, $\sqrt[N]{2^3}$, $\dots$, $\sqrt[N]{2^{N-1}}$ sono irrazionali, come stabilito dalla seguente

\noindent
{\bf Proposizione}. Se $m\geq 2$ e $n\geq 2$ sono numeri naturali e $m$ non \`e una potenza $n$--esima perfetta (ovvero $m$ non \`e del tipo $a^n$ per qualche numero $a$), allora $\sqrt[n]{m}$ \`e un numero irrazionale.

\noindent
{\bf Dim.} Supponiamo per assurdo che $\sqrt[n]{m}$ sia razionale: 
\begin{equation}
\label{ipass}
\sqrt[n]{m}=\dfrac{p}{q}, \quad p,q\;\;\;\textrm{primi tra loro}
\end{equation}
e scriviamo la scomposizione in fattori primi di $m$, $p$ e $q$:
$$
m=m_1^{k_1}m_2^{k_2}\cdot\cdot\cdot m_r^{k_r}, \quad 
p=p_1^{h_1}p_2^{h_2}\cdot\cdot\cdot p_s^{h_s}, \quad 
q=q_1^{l_1}q_2^{l_2}\cdot\cdot\cdot q_t^{l_t}.
$$
Elevando ad esponente $n$ la (\ref{ipass}), si ha $mp^n=q^n$ e sostituendo 
le scomposizioni in fattori si trova
\begin{equation}
\label{scomp}
m_1^{k_1}m_2^{k_2}\cdot\cdot\cdot m_r^{k_r} 
p_1^{nh_1}p_2^{nh_2}\cdot\cdot\cdot p_s^{nh_s}=  
q_1^{nl_1}q_2^{nl_2}\cdot\cdot\cdot q_t^{nl_t}
\end{equation}
In base alle premesse:
\begin{itemize}
\item[$(i)$]  almeno uno degli esponenti $k_j$ tra $k_1$, $\dots$, $k_r$ deve essere non divisibile per $n$,
dato che $m$ non \`e una potenza ennesima perfetta,
\item[$(ii)$] nessuno dei fattori di $p_1$, $p_2$, $\dots$, $p_s$ pu\`o coincidere con uno dei fattori $q_1$, $q_2$, $\dots$, $q_t$, dato che $p$ e $q$ sono primi tra loro,
\end{itemize}
ragioniamo sull'espressione (\ref{scomp}). Trattandosi di fattori primi, i fattori $m_1$, $\dots$, $m_r$ devono essere tutti presenti a destra dell'uguale, ovvero comparire tra i fattori $q_1$, $\dots$, $q_t$, pertanto i fattori di $p$ sono tutti distinti dai fattori di $m$, per la premessa $(ii)$. 
Siamo giunti all'assurdo, in quanto per la premessa $(i)$ il fattore $m_j$ ha esponente $k_j$ non divisibile per $n$, mentre a destra dell'uguale tutti gli esponenti dei fattori sono multipli di $n$.
$\quad\square$

\noindent
Applicando il risultato della Proposizione ai numeri $\sqrt[N]{2}$, $\sqrt[N]{2^2}$, $\sqrt[N]{2^3}$, $\dots$, $\sqrt[N]{2^{N-1}}$ per $N>3$, poniamo $n=N$ e $m=2, 2^2, 2^3, \dots, 2^{N-1}$: dato che nessuno dei valori di $m$ pu\`o essere una potenza $N$--esima perfetta, si ha che tali coefficienti sono tutti \emph{numeri irrazionali}.

\subsection{Propriet\`a della scala equabile}

\noindent
I suoni della scala \fbox{E}, essendo equidistanti, presentano una caratteristica esclusiva: 

\noindent{\bf Propriet\`a di trasposizione degli indici}. L'insieme di suoni $\{f_{i_1}, \;\dots, \;f_{i_n}\}$ appartenenti alla scala equabile \`e congruente ad ogni 
insieme ottenuto sommando o sottraendo il me\-de\-si\-mo numero $k$ ad ogni indice:
\begin{equation}
\label{traspind}
\left\{ f_{i_1+k}, f_{i_2+k}, \dots f_{i_n+k}\right\} \equiv \left\{ f_{i_1}, f_{i_2}, \dots, f_{i_n}\right\}.
\end{equation}
Il riscontro \`e immediato: 
va verificata la (\ref{congruenti}) ponendo
$\left\{ {\widetilde f}_{i_1}, \dots, {\widetilde f}_{i_n}\right\}=
\left\{ f_{i_1+k}, \dots f_{i_n+k}\right\}$ e utilizzata la (\ref{ognifreq}), per scrivere
$\dfrac{f_{i_1+k}}{f_{i_1}}=\dfrac{r^{i_1+k}f_0}{r^{i_1}f_0}=r^k$, 
$\dfrac{f_{i_2+k}}{f_{i_2}}=\dfrac{r^{i_2+k}f_0}{r^{i_2}f_0}=r^k$ 
e cos\`i$\,$ via per ogni indice.

\noindent
Per andare alla ricerca di una lettura in termini musicali della propriet\`a di trasposizione \`e necessario 
precisare alcune nozioni di carattere musicale. 

\noindent
La definizione di intervallo (\ref{interv}) esige la presenza di due suoni, che possono essere intesi in successione (\emph{intervallo melodico}) oppure simultanei (\emph{intervallo armonico}). 
In senso pi\`u generale, l'insieme (\ref{seq}) pu\`o rappresentare una sequenza temporale di suoni 
(\emph{melodia}) oppure una sovrapposizione di suoni:
in quest'ultimo caso, due suoni sovrapposti formano una \emph{diade} o  \emph{bicordo}, per $n=3$ si hanno le \emph{triadi}, per $n=4$ le \emph{quadriadi}, queste ultime di largo impiego nella musica \emph{jazz}.
Per \emph{accordo} si intende una sovrapposizione di almeno tre suoni. In modo conforme alla consuetudine musicale, esprimiamo un bicordo o un accordo mediante gli intervalli che lo formano: per esempio la triade ${\cal A}=\{f_5, f_7, f_{11}\}$ disposta sulla scala $N=12$ \`e formata dai due intervalli $I_1=I(f_5, f_7)$ e $I_2=I(f_7, f_{11})$. In generale l'accordo di $n$ suoni $\{f_{i_1}, f_{i_2}, \dots, f_{i_n}\}$, $n\geq 3$, \`e formato e indicato dagli $n-1$ intervalli ordinati:
\begin{equation}
\label{intervacc}
I_1=I(f_{i_1}, f_{i_2}), \quad I_2=(f_{i_2}, f_{i_3}), \quad \dots \quad I_{n-1}=I(f_{i_{n-1}}, f_{i_n})
\end{equation}

\noindent
In tale contesto risulta chiaro, in base alla definizione (\ref{congruenti2}), che due melodie con egual numero di suoni $i_n$ sono congruenti se percorrono, nell'ordine, i medesimi intervalli. Allo stesso modo due bicordi o due accordi formati dal me\-de\-si\-mo numero di suoni sono congruenti se risultano formati dai medesimi intervalli (\ref{intervacc}). 

\noindent
La propriet\`a (\ref{traspind}) di trasporre gli indici sulla scala \fbox{E}, letta a vari livelli, comporta pi\`u conseguenze: 
\begin{itemize} 
\item[I.] trasponendo tutti gli indici della scala equabile $(f_0, 2f_0)$ di una medesima quantit\`a $K$ si ottiene la scala $\{f_K, \;f_{K+1}, \dots, f_{2K}\}$ ancora equabile;
\item[II.] la melodia $\{ f_{i_1\pm k}, \dots, f_{i_n \pm k}\}$  ottenuta trasponendo $\{f_{i_1}, \dots f_{i_n}\}$ \`e congruente alla melodia di partenza, qualunque sia $k$; 
\item[III.] A meno di congruenze, gli accordi con i me\-de\-si\-mi valori (\ref{intervacc}) sono identici.
\end{itemize}

\noindent
Dal punto di vista della composizione musicale, il suono $f_0$ d\`a la cosiddetta \emph{tonalit\`a} al brano: dunque, a meno di trasportare tutti i suoni pi\`u in alto o pi\`u in basso, secondo il medesimo termine di trasposizione, la scala equabile propone una sola tonalit\`a, a meno di congruenze.

\noindent
Analogamente, una melodia trasposta, ovvero eseguita in due tonalit\`a differenti, percorre intervalli di medesima lunghezza, ovvero non subisce alcuna alterazione dal punto di vista della struttura, ma solo dal punto di vista dell'altezza dei suoni.
L'omogeneit\`a della scala equabile si presenta anche dal punto di vista armonico, degli accordi, che trasportati in altra tonalit\`a non mutano gli elementi distintivi.

\noindent
Il commento finale rileva due vantaggi dal punto di vista delle regole, essendo agevole il fatto di 
 comporre in ogni tonalit\`a, senza problematiche specifiche per qualcuna di esse, 
e di trasportare suoni, melodie e armonie nelle tonalit\`a desiderate, con semplici regole aritmetiche.
Si pu\`o tuttavia cogliere uno svantaggio dal punto di vista del colore, della variet\`a del suono
dal momento che le tonalit\`a sono tutte congruenti, ovvero tutte uguali a meno di trasposizioni.

\noindent
Dal punto di vista della terminologia musicale l'omogeneit\`a scala equabile 
permette una collocazione decisamente semplice dei tipi di intervallo, di scala o di accordo: 
definiamo \emph{semitono} l'intervallo $I_{\cal S}$ tra due suoni con indici consecutivi 
$I_{\cal S}=\dfrac{f_{i+1}}{f_i}$
dunque coincidente con un gradino (\ref{gradi}), 
\emph{tono} l'intervallo $I_{\cal T}$ composto da due semitoni consecutivi, nel senso di (\ref{sommainterv}):
${\cal I}_{\cal T}={\cal I}_{{\cal S}_1}{\cal I}_{{\cal S}_2}=\dfrac{f_{i+1}}{f_i}\dfrac{f_{i+2}}{f_{i+1}}=\dfrac{f_{i+2}}{f_i}$. Poniamo $N=12$, coerentemente con il sistema musicale odierno: la scala dei semitoni ${\cal S}_C=\{ f_0, f_1, \dots, f_{11}, 2f_0\}$, qualunque sia $f_0$ \`e la \emph{scala cromatica}, il sottoinsieme ${\cal S}_D=\{ f_0, f_2, f_4, f_5, f_7, f_9, f_{11}, 2f_0\}$, in cui si alternano toni e semitoni, \`e la \emph{scala diatonica maggiore}. 
In quest'ultima il suono iniziale $f_0$ d\`a il nome alla scala medesima: ad esempio, se $f_0$ \`e un DO, allora la sequenza \`e la {\it scala di DO maggiore}. Gli altri termini abbinati ai gradi della scala diatonica sono sicuramente noti: ${\cal S}_D=\{ f_0, f_2, f_4, f_5, f_7, f_9, f_{11}, 2f_0\}=\{ DO, RE, MI, FA, SOL, LA, SI, DO\}$. L'aggiunta delle operazioni $\sharp$ (\emph{diesis}) e $\flat$ (\emph{bemolle}) 
di traslazione in alto e in basso, rispettivamente, di un semitono:
$$
\sharp (f_i)=f_{i+1}, \quad \flat(f_i)=f_{i-1}
$$
permette di nominare i gradi della scala cromatica: 
${\cal S}_C=\{DO, DO\sharp= RE\flat, RE, RE\sharp=MI\flat, MI, FA,$ $FA\sharp=SOL\flat, SOL, SOL\sharp=LA\flat,  LA, LA\sharp=SI\flat, SI, DO\}$.

\noindent
La classificazione degli intervalli, in base a quanto stabilito dalla propriet\`a (\ref{traspind}), pu\`o fare riferimento unicamente al numero di semitoni compresi tra gli estremi: limitiamoci per il momento a segnalare l'intervallo di \emph{quinta} $I(f_k, f_{k+7})$ con $7$ semitoni e l'intervallo di \emph{quarta} $I(f_k, f_{k+5})$, $5$ semitoni, qualunque sia $k$ intero. Ad esempio $I(f_0, f_7)=(DO, SOL)$, $I(f_3, f_{10})=(MI\flat, SI\flat)$, 
$I(f_6, f_{13})=(FA\sharp, DO\sharp)$ sono intervalli di quinta, tutti congruenti; allo stesso modo $I(f_0, f_4)=(DO, FA)$, $I(f_{11}, f_{15})=(LA\sharp, RE\sharp)$ sono intervalli di quarta.

\noindent
Indichiamo anche alcuni accordi: $\{f_k, f_{k+4}, f_{k+7}\}$ \`e la \emph{triade o accordo maggiore}, $\{ f_k, f_{k+3}, f_{k+7}\}$ \`e la \emph{triade o accordo minore}, $\{ f_k, f_{k+4}, f_{k+7}, f_{k+11}\}$ \`e la 
\emph{quadriade o accordo di settima maggiore}. La nota pi\`u grave d\`a il nome all'accordo: ad esempio, $(DO,MI,SOL)$ \`e l'accordo di DO maggiore, $(RE,FA,LA)$ quello di RE minore, $(MI, SOL\sharp, SI, RE\sharp)$ l'accordo di MI settima maggiore.

\subsection{Una legge sperimentale a supporto della scala equabile}

\noindent
La scala equabile, che abbiamo presentato per prima, data la sua semplicit\`a strutturale, \`e l'ultima in ordine di tempo tra i sistemi di suoni fondamentali proposti sin dall'antichit\`a: come vedremo pi\`u avanti, le scale precedenti si basano sulla commensurabilit\`a, ovvero sui numeri razionali.

\noindent
Un aspetto contestato alla scala \fbox{E} \`e stato proprio quello di imporre l'uso di numeri complicati dal punto di vista aritmetico: i numeri irrazionali hanno uno sviluppo decimale infinito, non periodico, necessitano di essere approssimati, ...  Pertanto la scala non fu unanimamente associata ad una naturale e spontanea condizione, come poteva essere quella di ottenere i suoni per divisioni in poche parti (dunque facendo uso di semplici numeri razionali che esprimono le divisioni) di una corda.

\noindent
In secondo luogo, se abbiamo presente la formazione fisica di un suono come sovrapposizione di suoni armonici, ovvero di suoni con frequenze $f$, $2f$, $3f$, $\dots$, \`e da notare l'estraneit\`a dei suoni della scala \fbox{E} dal punto di vista degli armonici: nessun suono, a parte ovviamente l'ottava, fa parte degli armonici di qualche altro.

\noindent
Tuttavia, in ambito sperimentale esiste una legge empirica (formulata gi\`a dal 1860) che si adatta perfettamente a spiegare 
la costruzione della scala equabile: tale legge, nota come {\it legge di Weber}, riguarda in generale la relazione tra uno stimolo (che pu\`o essere un peso da sopportare, un agente che provoca dolore, oppure, appunto, una fonte sonora) e la percezione che consegue. 
Si tratta evidentemente di rappresentare un fenomeno a carattere soggettivo (si parla infatti di psicofisica): in modo generale, sulla base di test ed esperimenti, tale fenomeno viene inquadrato da una formula matematica. 

\noindent
Se chiamiamo $S$ l'intensit\`a dello stimolo e $P$ la percezione dovuta allo stimolo, la legge afferma che
\begin{equation}
\label{weber}
S_1\Delta P=k\Delta S
\end{equation}
dove: $\Delta S=S_2-S_1$ \`e la variazione dello stimolo $S$, mentre $\Delta P=P_2-P_1$ \`e la variazione nella percezione, $k$ \`e una costante legata al contesto sperimentale. Per esempio, se si tratta di valutare la percezione di fatica ($P$) dovuta ad un peso da sorreggere ($S$), $\Delta S$ \`e il peso aggiunto a $S_1$ per arrivare a $S_2$, in modo da poter giudicare la nuova impressione percettiva $P_2$, quando invece $P_1$ \`e quella corrispondente allo scenario iniziale $S_1$.
In sostanza la legge afferma che la percezione non \`e graduale, uniforme rispetto allo stimolo, ma dipende 
dallo stato $S_1$ da cui si parte; questo fatto \`e comprensibile se rimaniamo nel contesto esemplificativo  del peso: a parit\`a di $\Delta S$, ovvero ``medesima aggiunta sul carico iniziale $S_1$'', se quest'ultimo \`e ingente la risposta percettiva $\Delta P$ (che riporta l'aumento di fatica) \`e inferiore, rispetto al caso in cui $S_1$ \`e lieve.

\noindent
Vogliamo dedurre da questa legge la seguente

\noindent
{\bf Propriet\`a}
L'aumento uniforme della percezione avviene in corrispondenza di stimoli che formano una progressione geometrica.

\noindent
{\bf Dim.} Partiamo da (\ref{weber}): aumento uniforme dei $P_1$, $P_2$, $P_3$, $\dots$, $P_N$ significa  $\Delta P$ tutti uguali: 
$$
\Delta P=P_2-P_1=P_3-P_2=\dots P_N-P_{N-1}=C,
$$
quest'ultimo valore costante. Sostituendo in (\ref{weber}), si trova 
$$
CS_j=k(S_{j+1}-S_j), \quad j=1,2,\dots N-1
$$
e, dividendo per $kS_j$: $\dfrac{S_{j+1}}{S_j}=1+\dfrac{C}{k}$, quest'ultimo valore costante, dunque la progressione $S_1$, $S_2$, $\dots$, $S_N$ \`e geometrica, di ragione $1+\dfrac{C}{k}$. $\quad\square$

\noindent
La Propriet\`a appena dimostrata fornisce una motivazione, dal punto di vista dello studio delle percezioni, alla caratteristica equidistanza della scala \fbox{E}: se infatti $S$ \`e lo stimolo uditivo ad una determinata frequenza $f$ (ora stiamo valutando la risposta all'altezza del suono, ma la medesima legge ha una pertinente applicazione anche all'intensit\`a del suono) e $P$ la risposta uditiva, la legge afferma che un graduale aumento della sensazione in altezza avviene in corrispondenza di frequenze che si susseguono in progressione geometrica, proprio come nella scala \fbox{E}. 
Pur rimanendo nella sfera delle percezioni e non delle leggi fisiche automaticamente quantificabili, l'obiettivo di avvertire un aumento dell'altezza dei suoni progressivo ed uniforme viene realizzato dalla scala dei suoni irrazionali.

 
\section{La scala razionale che permette due soli divisori: scala pitagorica}

\subsection{Il Principio della quinta}

\noindent
Il clima \`e ora completamente diverso: alla praticit\`a strutturale di separare in modo equidistante (nel senso del rapporto) gli elementi della scala, si sostituisce un principio che sostanzialmente estende quello dell'ottava: quest'ultimo, ricordando la (\ref{2kf}), impone la necessit\`a di avere sulla scala i suoni ottenuti per divisione a met\`a della corda, o per raddoppio della corda. 
L'estensione del Principio accoglie la presenza, nella scala da costruire, del suono che si ottiene dividendo in tre parti uguali la corda o triplicandone la lunghezza.
Enunciamo questa regola attraverso il 

\noindent
{\bf Principio della quinta}.  Possono entrare a far parte della scala che contiene il suono $f$ i suoni 
$3f$, $9f$, $27f$, $3^Nf$, $\dots$ e i suoni 
$\frac{1}{3}f$, $\frac{1}{9} f$, $\frac{1}{27}f$, $\dots$, $\frac{1}{3^N}$, $\dots$
per ogni $N=1,2,3,\dots$. 
Ovvero, la presenza di $f$ nella scala permette di includere nella medesima scala i suoni 
\begin{equation}
\label{3kf}
3^kf, \quad k\in {\Bbb Z}\quad \textrm{numeri interi}.
\end{equation} 

\noindent
Dal punto di vista delle motivazioni, possiamo senz'altro ripercorrere alcune di quelle riportate dopo la (\ref{2kf}): il suono prodotto dalla corda divisa in tre risulta particolarmente consonante con quello generato dalla corda intera: subito dopo l'ottava, nella gerarchia delle consonanze va collocata questa circostanza; d'altra parte, in un senso filosofico, se \`e il Numero a regolare ogni aspetto dell'Universo, dunque anche nella musica il numero tre, dopo il due, prescrive certamente un dettame primario.  
Da un punto di vista invece fisico, nella successione delle armoniche, il modo naturale di vibrazione della corda in tre parti uguali segue immedatamente quello della vibrazione in due parti.

\noindent
\`E opportuno chiarire la denominazione \emph{quinta}: applichiamo il Principio della quinta ad un suono $f_0$: il suono $3f_0$ deve essere incluso e, riportato sull'ottava $[f_0, 2f_0]$ corrisponde a $\dfrac{3}{2}f_0=1,5f_0$. Operando un confronto con i suoni scala equabile (\ref{ognifreq}) riscontriamo che $\sqrt[12]{2^7}f_0\approx 1,49830$ \`e il suono maggiormente vicino a $\dfrac{3}{2}f_0$ e lo scarto \`e poco pi\`u che percettibile.
La lunghezza pressoch\'e identica di $I(f_0, f_7)$, $f_7$ suono equabile, e di $I\left(f_0, \dfrac{3}{2}f_0\right)$  fa utilizzare lo stesso termine per definire l'intervallo.
Definiamo quindi 
\emph{quinta pitagorica} ogni intervallo $I$ per cui 
\begin{equation}
\label{quintapit}
I(f_{i_1}, f_{i_2})=\dfrac{f_{i_2}}{f_{i_1}}=\dfrac{3}{2}
\end{equation}
Calcoliamo l'ampiezza dell'intervallo complementare alla quinta nell'arco dell'ottava:
$I\left(\frac{3}{2}f_0, 2f_0\right)=\frac{4}{3}$. 
Passando in rassegna i valori (\ref{ognifreq}), si vede che $\frac{4}{3}=1,333...$ ha come valore pi\`u vicino $f_5=1,33484...$
La denominazione \emph{intervallo di quarta} gi\`a assegnata a $(f_0, f_5)$ fa comprendere la definizione di 
\emph{quarta pitagorica}, come ogni intervallo $I$ per cui 
\begin{equation}
\label{quartapit}
I(f_{i_1}, f_{i_2})=\dfrac{f_{i_2}}{f_{i_1}}=\dfrac{4}{3}
\end{equation}
Ci siamo appena scontrati con un fatto consueto nella comparazione di pi\`u tipi di scala: un intervallo, pur qualificato con lo stesso nome, pu\`o avere ampiezza diversa, se collocato su scale di diversa formazione. 
In modo preciso, dunque, dovremo parlare di \emph{quinta, quarta equabile}, oppure di \emph{quinta, quarta pitagorica}, a seconda della scala.

\noindent
Storicamente, la scala costruita mediante (\ref{3kf}) va senz'altro ricondotta alle teorie pitagoriche sulla musica. Il nome assegnato all'intervallo $I\left(f_0, \frac{3}{2}f_0\right)$ era quello di \emph{diapente}, ovvero ``attravero cinque toni''. L'intervallo complementare $I\left(\frac{3}{2}f_0, 2f_0\right)$ corrisponde  alla quarta, denominata '\emph{diatessaron}, ovvero ``attraverso quattro toni''. 
L'operazione di formazione dell'ottava da una quinta ed una quarta: \emph{diapason}=\emph{diapente}+ \emph{diatessaron}, che noi traduciamo nei termini di (\ref{sommainterv}) come
$I(f_0, 2f_0) = I\left(f_0, \frac{3}{2}f_0\right) I\left(\frac{3}{2}f_0, 2f_0\right)$ e che rappresentiamo graficamente come {\footnotesize $
\begin{array}{l}
\overbrace{
\overbrace{|-----|}^{I=\frac{3}{2}}-----|}^{I=2}\\
f_0 \qquad \quad \quad  \underbrace{\frac{3}{2}f_0 \;\quad \quad \quad 2f_0}_{I=\frac{4}{3}}
\end{array}
$
} ha la sua antica raffigurazione tramite semicerchi:

{\footnotesize 

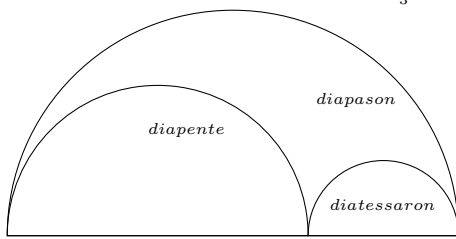
\begin{tikzpicture}
\draw[] (0,0) arc (180:0:3) -- cycle;
\draw[] (0,0) arc (180:0:2) -- cycle;
\draw[] (4,0) arc (180:0:1) -- cycle;
\coordinate [label=below left:\textcolor{black}{$\scriptstyle diapente$}](a) at (3,1.6);
\coordinate [label=below left:\textcolor{black}{$\scriptstyle diatessaron$}](a) at (5.8,.6);
\coordinate [label=below right:\textcolor{black}{$\scriptstyle diapason$}](b) at (4,2);
\coordinate (c) at ($(a)!{1/6}!(b)$);;
\end{tikzpicture}
}

\subsection{La scala pitagorica}

\noindent
Operiamo ora la costruzione della scala basata sulla regola (\ref{3kf}), oltre che, naturalmente, sulla (\ref{2kf}). 
Pertiamo dal suono $f_0$: 
evidentemente, le divisioni o moltiplicazioni per $2$ o per $3$ delle due regole portano prendere in esame i suoni
\begin{equation}
\label{scalapit}
f=2^{m} 3^{n}f_0, \qquad m, n \in {\Bbb Z}
\end{equation}
dove i numeri interi $m$ e $n$ possono essere positivi o negativi. 
In particolare, rendiamo esplicite le seguenti operazioni, la prima due delle quali \`e gi\`a nota:
\begin{itemize}
\item[$(i)$] moltiplicare [risp.~divedere] per $2$ significa portare all'ottava sopra [sotto]: l'esponente positivo [negativo] $m$ al fattore $2$ fa ripetere l'operazione, ovvero porta il suono $m$ ottave in avanti, verso l'acuto [$-m$ ottave all'indietro, verso il grave];

\item[$(ii)$] ogni moltiplicazione per $3$ porta il suono alla quinta dell'ottava sopra, ovvero lo traferisce un'ottava ed una quinta sopra: 

\item[$(iii)$] ogni divisione per $3$ porta il suono un'ottava ed una quinta sotto:
\end{itemize}
Le operazioni $(ii)$ e $(iii)$ possono essere schematizzate come segue: nella prima, a sinistra, si passa da $f_0$ a $3f_0$, nella seconda da $f_0$ a $\frac{1}{3}f_0$. Si evidenziano anche gli intervalli di quarta che 
completano le ottave in cui si estendono le operazioni.

{\footnotesize
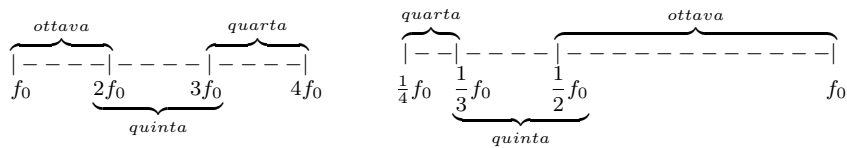
$$
\begin{array}{l}
\overbrace{|----|}^{ottava}----\overbrace{|----|}^{quarta} \\
f_0 \;\;\qquad \underbrace{2f_0\;\; \qquad \;3f_0}_{quinta} \;\;\;\qquad 4f_0 
\end{array}
\qquad
\begin{array}{l}
\; \overbrace{|--|}^{quarta}----\overbrace{|------------|}^{ottava}\\
\frac{1}{4}f_0\; \;\,\underbrace{\frac{1}{3}f_0 \;\;\quad \;\; \,\;\frac{1}{2}f_0}_{quinta}  \hspace{3.1truecm} f_0
\end{array}
$$
}
L'esponente $n$ in (\ref{scalapit}) itera l'operazione $(ii)$ se positivo, oppure l'operazione $(iii)$, se negativo.

\subsection{Le propriet\`a dei numeri della scala pitagorica}

\noindent
A partire da un  suono $f_0$, nel caso della scala equabile i suoni (\ref{2kf}) sono automaticamente incorporati nella scala, mentre i suoni intermedi provengono dall'ipotesi di equidistanza. Nel caso in esame, invece, la (\ref{2kf}) e la (\ref{3kf}) ammette gi\`a un grande numero di suoni e non \`e necessario farne intervenire altri nella scala: \`e dunque opportuno esaminare le propriet\`a dei numeri (\ref{scalapit}).

\noindent
{\bf Proposizione}.  I numeri (\ref{scalapit})
\begin{itemize}
\item[$(i)$]  sono razionali, 
\item[$(ii)$] sono infiniti e tutti distinti tra loro, se gli esponenti sono distinti,
\end{itemize}

\noindent
{\bf Dim}.  La propriet\`a $(i)$ \`e evidente: si tratta di numeri ottenuti mediante prodotti o divisioni con numeri interi. Per la $(ii)$: se ipotizziamo
$2^{m_1}3^{n_1}=2^{m_2}3^{n_2}$ con $m_1, m_2, n_1, n_2 \in {\Bbb Z}$, troviamo 
\begin{equation}
\label{23}
2^{m_1-m_2}=3^{n_2-n_1}.
\end{equation}
Vogliamo dimostrare che la (\ref{23}) \`e compatibile solo se $m_1=m_2$ e $n_1=n_2$. Infatti, se fosse $m_1>m_2$, anche 
$n_2$ deve essere maggiore di $n_1$, dato che $2^{m_1-m_2}>1$. Tuttavia, la (\ref{23}) porrebbe uguali un numero dispari ed un numero pari. D'altra parte, se $m_1<m_2$, si conclude allo stesso modo, ragionando  su gli inversi  $2^{m_2-m_1}$ e $3^{n_2-n_1}$. Dunque i numeri di tipo (\ref{scalapit}) sono uguali se e solo se gli esponenti sono ordinatamente uguali. $\quad\square$

\noindent
Come caso particolare della propriet\`a $(ii)$, enunciamo il 

\noindent
{\bf Corollario}. Non \`e possibile determinare $m\not= 1$ e $n\not= 0$ in modo che (\ref{scalapit}) dia il numero $2$.

\noindent
Il Corollario va ricondotto alla situazione nota come ``ciclo delle quinte che non si  chiude'': se interviene almeno una quinta ($n\not =2$), non si potr\`a mai tornare su un'ottava del suono di base $f_0$.

\noindent
Se chiamiamo \emph{scala pitagorica} o \emph{scala} \fbox{P} l'insieme dei suoni (\ref{scalapit}), la principale caratteristica della scala \`e riassunta dalla seguente

\noindent
{\bf Propriet\`a}. La scala \fbox{P} \`e una scala infinita di suoni tutti diversi.

\subsection{Un numero finito di suoni per l'ottava della scala pitagorica}

\noindent
La propriet\`a matematica dei numeri generati dal $2$ e dal $3$ tramite (\ref{scalapit}) si traduce nella 
ricchezza illimitata di suoni a disposizione.
Se, dal punto di vista della musica teorica,  questa connotazione pu\`o apparire avvincente, dal punto di vista della musica pratica, soprattutto strumentale, si delineano degli inconvenienti pesanti, subiti ad esempio dagli strumenti a tastiera, sui quali non esiste la possibilit\`a di produrre tutte le frequenze possibili.
Per porre un termine all'infinita produzione di suoni da parte di (\ref{scalapit}), va escogitato un criterio iterativo, da interrompere quando si ritene opportuno.

\noindent
In questo compito di delineare una modalit\`a di selezione di un numero finito di suoni (\ref{scalapit}), siamo avvantaggiati dalla conoscenza della scala dei $12$ suoni equabili, con i quali possiamo operare un confronto, come \`e gi\`a  avvenuto per l'intervallo di quinta e di quarta.
Tuttavia, non \`e ora scontato che $N$ (numero dei suoni nell'ottava) sia $12$: anzi, procedere un p\`o oltre il $12$ dar\`a luogo alla maggiore ricchezza di suoni che la voce umana e alcuni strumenti riescono a produrre.
La procedura di selezione per \fbox{P} \`e senz'altro intuitiva: a partire da $f_0$, si selezionano i suoni tra i (\ref{scalapit})) avanzando e retrocedendo per intervalli di quinta pitagorica e riportando eventualmente i suoni nell'ottava di riferimento, se questa viene oltrepassata.
Con questa premessa, possiamo passare alla formulazione del procedimento: 

\begin{itemize}
\item[$I.$]
si considera la successione di suoni a sinistra e a destra di $f_0$ procedendo per intervalli di quinta pitagorica: 
\begin{equation*}
\left(\frac{3}{2}\right)^{-m_1} 
\dots \left(\frac{3}{2}\right)^{-3}f_0,
\;\left(\frac{3}{2}\right)^{-2}f_0, \;\;\left(\frac{3}{2}\right)^{-1}f_0 \;\leftarrow\;\textrm{\fbox{$f_0$}} \;
\rightarrow\;\dfrac{3}{2}f_0, \;\;
\left(\frac{3}{2}\right)^2f_0, \;\;\left(\frac{3}{2}\right)^3f_0 \; \dots \; \left(\frac{3}{2}\right)^{m_2}
\end{equation*}

\item[$II.$] si riportano  i suoni nell'ottava di riferimento $(f_0, 2f_0)$ mediante moltiplicazioni con le potenze del $2$.
\end{itemize}

\noindent
Operativamente, si tratta di effettuare i calcoli
\begin{equation}
\label{perquinte}
\begin{array}{l}
\textrm{per}\;k_1=1, \dots, m_1:\;
f_{h_1,k_1}=2^{h_1}\left(\frac{3}{2}\right)^{-k_1}f_0 \;\textrm{in modo che}\;
1<\frac{f_{h_1,k_1}}{f_0}<2,\;\;h_1\geq 0,\\[11pt]
\textrm{per}\;k_2=1,\dots, m_2:\;
f_{h_2,k_2}=2^{h_2}\left(\frac{3}{2}\right)^{k_2}f_0\;\textrm{in modo che}\;
1<\frac{f_{h_2,k_2}}{f_0}<2, \;\;h_2\leq 0.
\end{array}
\end{equation}

\noindent
Riportiamo nella seguente tabella il calcolo di alcuni valori (normalizzati con $f_0$) dopo $12$ quinte in avanti ($m_1=12$) e $12$ all'indietro ($m_2=12$), per ottenere quindi la selezione finita della scala pitagorica con $26$ suoni, se comprendiamo $f_0$ e $2f_0$:

\begin{center}
\begin{tabular}{|c||c|}\hline
 $f_{h_1,k_1}/f_0=2^{h_1}(\frac{3}{2})^{k_1}$ & $f_{h_2,k_2}/f_0=2^{h_2}(\frac{3}{2})^{k_2}$\\ \hline \hline
$2\left(\frac{3}{2}\right)^{-1}=\frac{4}{3}=1,{\bar 3}$ & $\frac{3}{2}=1,5$\\ \hline
 $2^2\left(\frac{3}{2}\right)^{-2}=\frac{16}{9}=1,{\bar 7}$ & $2^{-1}\left(\frac{3}{2}\right)^{2}=\frac{9}{8}=1,125$ \\ \hline
 $2^2\left(\frac{3}{2}\right)^{-3}=\frac{32}{27}=1,\overline{185}$ & $2^{-1}\left(\frac{3}{2}\right)^{3}=\frac{27}{16}=1,6875$ \\ \hline
 $2^3\left(\frac{3}{2}\right)^{-4}=\frac{128}{81}=1,58024...$ & $2^{-2}\left(\frac{3}{2}\right)^{4}=\frac{81}{64}=1,26562...$ \\ \hline
 $2^3\left(\frac{3}{2}\right)^{-5}=\frac{256}{243}=1,05349...$ & $2^{-2}\left(\frac{3}{2}\right)^{5}=\frac{243}{128}=1.89843...$ \\ \hline
 $2^4\left(\frac{3}{2}\right)^{-6}=\frac{1024}{729}=1,40466...$ & $2^{-3}\left(\frac{3}{2}\right)^{6}=\frac{729}{512}=1,42382...$ \\ \hline
 $2^5\left(\frac{3}{2}\right)^{-7}=\frac{4096}{2187}=1,87288...$ & $2^{-4}\left(\frac{3}{2}\right)^{7}=\frac{2187}{2048}=1,06787...$ \\ \hline 
 $2^5\left(\frac{3}{2}\right)^{-8}=\frac{8192}{6561}=1,24859...$ & $2^{-4}\left(\frac{3}{2}\right)^{8}=\frac{6561}{4096}=1,60180...$ \\ \hline
 $2^6\left(\frac{3}{2}\right)^{-9}=\frac{32768}{19683}=1,66478...$ & $2^{-5}\left(\frac{3}{2}\right)^9=\frac{19683}{16384}=1,20135...$ \\ \hline
 $2^6\left(\frac{3}{2}\right)^{-10}=\frac{65536}{59049}=1,10985...$ & $2^{-5}\left(\frac{3}{2}\right)^{10}=\frac{59049}{32768}=1,80203...$ \\ [2ex] \hline 
 $2^7\left(\frac{3}{2}\right)^{-11}=\frac{262144}{177147}=1,47981...$ & $2^{-6}\left(\frac{3}{2}\right)^{11}=\frac{177147}{131072}=1,35152...$ \\ [2ex] \hline 
 $2^8\left(\frac{3}{2}\right)^{-12}=\frac{1048576}{531441}=1,97308...$ & $2^{-7}\left(\frac{3}{2}\right)^{12}=\frac{531441}{524288}=1,01364...$ \\ [2ex] \hline 
\end{tabular}
\end{center}

\noindent
Da $k_1=k_2=4$ in poi abbiamo scritto solo le prime cinque cifre decimali, ma ricordiamo che si tratta in tutti i casi di numeri razionali (in particolare: quelli a sinistra sono periodici, quelli a destra hanno un numero finito di cifre decimali).

\noindent
La scelta di $12$ passi \`e giustificata dalla seguente osservazione: nell'ultimo passaggio $k_1=k_2=12$ i valori ottenuti sono molto prossimi a $1$ (se in avanti) e a $2$ (se all'indietro), lasciando pensare che ritroviamo suoni probabilmente confondibili con $f_0$ e $2f_0$, rispettivamente. Tale rispondenza offre la possibilit\`a di chiudere la scala, prevedendo che dal dodicesimo passo in poi i calcoli portino a valori 
quasi coincidenti con quelli gi\`a trovati.

\noindent
Dobbiamo ora occuparci di operare un confronto fra i valori calcolati e quelli della scala equabile (\ref{ognifreq}), con $N=12$: i ventiquattro suoni della tabella precedente, ai quali si aggiungono gli estremi $f_0$ e $2f_0$, vengono incasellati in corrispondenza del suono equabile pi\`u vicino (come nella tabella precedente, i valori sono normalizzati con $f_0$).

\begin{center}
\begin{tabular}{||l||l||}\hline
scala \fbox{E} & scala \fbox{P} \\ \hline
$*\;1$ & $1\qquad \quad  \quad \;\;\;2^{-7}(\frac{3}{2})^{12}$ \\ \hline 
$\;\;\sqrt[12]{2}=1,05946...$ &  $ 2^3(\frac{3}{2})^{-5} \qquad 2^{-4}(\frac{3}{2})^7$ \\ \hline
$*\;\sqrt[12]{2^2}=1,12246...$ & $2^{-1}(\frac{3}{2})^2\qquad 2^6(\frac{3}{2})^{-10}$ \\ \hline
$\;\;\sqrt[12]{2^3}=1,18921...$ & $2^2(\frac{3}{2})^{-3} \qquad 2^{-5}(\frac{3}{2})^9$ \\ \hline
$*\;\sqrt[12]{2^4}=1,25992...$ & $2^{-2}(\frac{3}{2})^4\qquad 2^5(\frac{3}{2})^{-8}$ \\ \hline 
$*\;\sqrt[12]{2^5}=1,33484...$ & $2(\frac{3}{2})^{-1}\qquad \;\;2^{-6}(\frac{3}{2})^{11}$ \\ \hline 
$\;\;\sqrt[12]{2^6}=1,41421...$ & $2^4(\frac{3}{2})^{-6} \qquad 2^{-3}(\frac{3}{2})^6$ \\ \hline 
$*\;\sqrt[12]{2^7}=1,49830...$ & $\frac{3}{2} \qquad \quad \quad \;\;\; 2^7(\frac{3}{2})^{-11}$ \\ \hline 
$\;\;\sqrt[12]{2^8}=1,58740...$ & $2^3(\frac{3}{2})^{-4} \qquad 2^{-4}(\frac{3}{2})^8$ \\ \hline 
$*\;\sqrt[12]{2^9}=1,68180...$ & $2^{-1}(\frac{3}{2})^3\qquad 2^6(\frac{3}{2})^{-9}$ \\ \hline 
$\;\;\sqrt[12]{2^{10}}=1,78180...$ & $2^2(\frac{3}{2})^{-2} \qquad 2^{-5}(\frac{3}{2})^{10}$ \\ \hline 
$*\;\sqrt[12]{2^{11}}=1,88775...$ & $2^{-2}(\frac{3}{2})^5  \qquad 2^5(\frac{3}{2})^{-7}$ \\  \hline 
$*\;2\;$ & $2 \qquad \quad \quad \;\; \;2^8(\frac{3}{2})^{-12}$ \\ \hline 
\end{tabular}
\end{center}
Con l'asterisco sono contrassegnati i suoni della scala diatonica. 
L'assetto dei suoni pitagorici attorno a quelli equabili si presenta ordinata e sistematica: 
i $26$ suoni si dispongono a coppie offrendo un'approssimazione per difetto ed una per eccesso di ciascuno dei $13$ suoni equabili. Nell'ordinare ciascuna coppia su ogni riga della tabella si \`e rispettata la seguente logica: per i suoni diatonici $*$ si \`e scritto per primo il valore incontrato con un minor numero di iterazioni (ovvero con $k_1$ o $k_2$ inferiore), mentre per le restanti coppie il primo elemento della coppia corrisponde al valore per difetto, il secondo al valore per eccesso.

\subsection{La scala musicale pitagorica}

\noindent
Dal punto di vista della matematica il procedimento non pu\`o che considerarsi soddisfacente: un numero stabilito di valori (\ref{scalapit}) riproduce sostanzialmente la scala equabile, in una versione duplicata di valori leggermente in difetto e valori leggermente in difetto.
Dal punto di vista della musica il sistema pitagorico non \`e tuttavia esattamente conformato all'esito della procedura, ma compie, per motivi pratici, una selezione sui $26$ suoni, ammettendo un solo valore per i suoni della scala diatonica, contrassegnati con $*$.
Per gli altri suoni, invece, si include sia il valore in eccesso, contrassegnato con $\sharp$, sia il valore in difetto, contrassegnato con $\flat$. I $18$ suoni 
pitagorici selezionati (compresi $f_0$ e $2f_0$) sono elencati nella seguente tabella,  evidenziandone la struttura (\ref{scalapit}) e abbinandoli ai corrispondenti nomi, se il suono di base \`e il DO:

\begin{center}
\begin{tabular}{||l|l||l | l ||}\hline
DO & $1$ & FA$\sharp$ & $2^{-9}3^6=1,42382...$ \\
RE$\flat$   &  $2^83^{-5}=1,05349...$ &  SOL   &  $ 2^{-1}3=1,5$\\
DO$\sharp$   &  $2^{-11}3^7=1,06787...$ &   LA$\flat$    &  $2^73^{-4}=1,58024...$\\
RE     &  $2^{-3}3^2=1,125$ &      SOL$\sharp$     &  $2^{-12}3^8=1,60180...$\\
MI$\flat$      &  $2^53^{-3}=1,\overline{185}$ &     LA      &  $2^{-4}3^3=1,6875$\\
RE$\sharp$         &  $2^{-14}3^9=1,20135...$ &    SI$\flat$       &  $2^43^{-2}=1,{\bar 7} $\\
MI   &  $2^{-6}3^4=1,26562...$ &     LA$\sharp$      &  $2^{-15}3^{10}=1,80203...$\\
FA    &  $2^23^{-1}=1,{\bar 3}$ &     SI      &  $ 2^{-7}3^5 = 1,89843...$\\
SOL$\flat$     &  $2^{10}3^{-6}=1,40466... $ &    DO       &  $2$\\
\hline
\end{tabular}
\end{center}
Si nota una differenza dei suoni bemolle e diesis interposti tra quelli diatonici.
Questa distinzione non esiste negli strumenti accordati mediante la scala equabile, come il pianoforte, dove non esiste la distinzione , per esempio, 
tra il do aumentato (DO$\sharp$) e il RE diminuito (RE$\flat$). 
Tuttavia, la maggiore ricchezza sonora della scala $\fbox{P}$ \`e riscontrabile in alcune circostanze: 
la voce umana oppure gli strumenti ad arco, come il violino, la viola, ..., in grado di produrre tutte le frequenze, possono manifestare l'attitudine a disporre alcuni intervalli secondo la scala pitagorica, come ad esempio l'intervallo di quinta, la cui misura pitagorica $3/2$ appare molto pi\`u naturale del valore leggermente inferiore e non razionale $\sqrt[12]{2^7}$. D'altra parte, la propensione ad eseguire in modo differente le note con alterazioni \`e un fatto naturale: il diesis $\sharp$, che tende a far proseguire verso l'acuto, \`e leggermente pi\`u alto del bemolle $\flat$, che tende a risolvere verso note pi\`u gravi. 
Un ascoltatore allenato pu\`o in effetti riscontrare, ad esempio per la voce o per un violino, la differenza fra due suoni che il sistema equabile dichiara identici.
Questa maggiore ricchezza rende senz'altro la scala pitagorica pi\`u flessibile alle sfumature 
del discorso musicale.

\noindent
Possiamo esaminare brevemente la situazione a livello di semitoni. 
La suddivisione del tono DO--RE, che nella scala equabile viene semplificata mediante un unico semitono, 
nella scala pitagorica \`e pi\`u articolata, come si pu\`o osservare dalla seguente schematizzazione:
$$
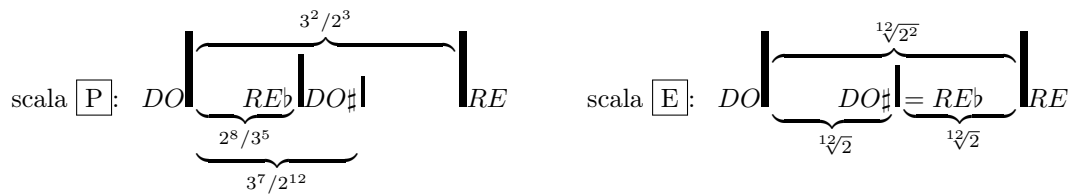
\begin{array}{ll}
\textrm{scala \fbox{P}:}\;\;\;
DO
\rule{1mm}{1cm}
\overbrace{
\underbrace{
\underbrace{
\hspace{.6cm}
RE\flat
}_{2^8/3^5}\rule{0.7mm}{.7cm}
DO\sharp
}_{3^7/2^{12}}\rule{0.4mm}{.4cm}
\hspace{1.2cm}
}^{3^2/2^3}\rule{1mm}{1cm}RE
&
\qquad \textrm{scala \fbox{E}:}\;\;\;
DO\rule{1mm}{1cm}\overbrace{
\underbrace{
\hspace{.86cm}
DO\sharp}_{\sqrt[12]{2}}\rule{0.55mm}{.55cm}
\underbrace{
=RE\flat
\hspace{.42cm}
}_{\sqrt[12]{2}}
}^{\sqrt[12]{2^2}}\rule{1mm}{1cm}RE
\end{array}
$$
Le disuguaglianze $\dfrac{2^8}{3^5}<\sqrt[12]{2}<\dfrac{3^7}{2^{12}}$ e $\sqrt[12]{4}<\dfrac{3^2}{2^3}$ fanno s\`\i$\,$che il semitono sulla scala equabile si posizioni all'interno dei due semitoni pitagorici DO--DO$\sharp$ e DO--RE$\flat$ e che il tono DO--RE sia pi\`u corto sulla scala equabile rispetto a quella  pitagorica.
Si pu\`o anche notare che gli intervalli RE$\flat$--RE e DO$\sharp$--RE sono esattamente la met\`a dei semitoni complementari dell'altro tipo: 
$I(RE\flat,RE)=\dfrac{3^2}{2^3}\dfrac{3^5}{2^8}=
\dfrac{3^7}{2^{11}}=2\times \dfrac{3^7}{2^{12}}=I(DO, RE\flat)$, 
$I(DO\sharp,RE)=\dfrac{3^2}{2^3}\dfrac{2^{12}}{3^7}=
\dfrac{2^9}{3^5}=
2\times \dfrac{2^8}{3^5}=I(DO,RE\flat)$.

\noindent
Dal punto di vista della realizzazione, 
uno svantaggio della scala pitagorica \`e la dipendenza di essa da $f_0$: 
la selezione della scala pitagorica attraverso un numero finito di suoni dipende dalla scelta di partenza $f_0$, ovvero un differente suono di base ${\widetilde f}_0$ d\`a origine, a parit\`a di passaggi in avanti e all'indietro, ad una scala differente. In effetti, se utilizziamo come suono di base ${\widetilde f}_0$ della scala di $f_0$, con  $f_0<{\widetilde f}_0<2f_0$ e procediamo, in base alla formulazione stabilita, $12$ quinte in avanti e $12$ all'indietro:
$$
{\widetilde f}_0,  \quad \left(\frac{3}{2}\right)^{\pm 1}{\widetilde f}_0, \quad \left(\frac{3}{2} \right)^{\pm 2}{\widetilde f}_0,\quad, \dots
$$
incontreremo sicuramente almeno un suono differente dalla scala di origine $f_0$, in virt\`u dell'unicit\`a della rappresentazione (\ref{scalapit}). Possiamo chiarire con un esempio: a partire da $f_0$ come nelle tabelle precedenti consideriamo ${\widetilde f}_0=2^{-6}3^4f_0$ (MI) ed utilizziamo quest'ultimo per dare origine al suono
${\widetilde f}_{h_2,k_2}=2^{-5}\left(\dfrac{3}{2}\right)^9{\widetilde f}_0$, $9$ quinte in avanti e $5$ suddivisioni a met\`a.
Tale suono non fa parte dei $26$ generati da 
$f_0=$DO: infatti
${\widetilde f}_{h_2,k_2}=2^{-20}3^{13}f_0=2^{-7}(\frac{3}{2})^{13}f_0$ ed occorrerebbero ulteriori iterazioni per raggiungerlo.

\noindent
Dunque, le scale pitagoriche finite non sono tutte congruenti fra loro, come avveniva per la scala equabile.
Solo se si considerassero tutti quanti i numeri (\ref{scalapit}), ovvero per gli infiniti valori degli esponenti, tutte le scale pitagoriche sarebbero uguali, indipendentemente dalla partenza.

\noindent
La mancanza di congruenza della scala pitagorica al variare del suono di base va inquadrata pi\`u in generale nella non validit\`a delle propriet\`a di trasposizione (\ref{traspind}): dal momento che le suddivisioni dell'ottava non sono equispaziate, non \`e detto che un'arbitraria traslazione degli indici trasferisca i suoni in altri gi\`a presenti sulla scala. Contrariamente a quanto si \`e rilevato mediante (\ref{traspind}) per la scala equabile, stavolta la trasposizione di una melodia non necessariamente fa percorrere i medesimi intervalli, un accordo trasportato non necessariamente \`e congruente con quello di partenza.
Dal punto di vista musicale, si afferma che le tonalit\`a non sono tutte uguali.

\section{La scala razionale con pi\`u divisori}

\noindent
Ripercorriamo ora i medesimi passi della costruzione della scala pitagorica, includendo un fattore primo in pi\`u nella categoria dei numeri (\ref{scalapit}): 
la base teorica aggiunge al Principio dell'ottava e della quinta il seguente

\noindent
{\bf Principio della terza}: la scala che contiene il suono $f$ pu\`o includere i suoni 
$5f$, $25f$, $125f$, $\dots$, $5^Nf$, $\dots$ e  
$\frac{1}{5}f$, $\frac{1}{25} f$, $\frac{1}{125}f$, $\dots$, $\frac{1}{3^N}$, $\dots$
per ogni $N=1,2,3,\dots$. In sintesi, nella scala possono essere presenti i suoni 
\begin{equation}
\label{5kf}
5^kf, \quad k\in {\Bbb Z}, \quad {\Bbb Z}\;\;\textrm{numeri interi}.
\end{equation}

\noindent
Di nuovo, una motivazione da accludere \`e il modo naturale di vibrazione di una corda in cinque parti uguali, successivo a tre e quattro, quest'ultimo gi\`a incluso nel Principio dell'ottava. Inoltre, ancora una volta, il termine \emph{terza} fa riferimento ai gradi della scala diatonica: consideriamo il suono $5f_0$ e riportiamolo nell'ottava $[f_0, 2f_0]$, per ottenere $\frac{5}{4}f_0=1,25f_0$. Percorrendo le tabelle precedenti, troviamo tale suono prossimo a $\sqrt[12]{2^4}f_0=1,25992... f_0$ sulla scala \fbox{E} e a $\dfrac{81}{64}f_0=1,26562... f_0$ sulla scala \fbox{P}, corrispondenti rispettivamente a $MI_{E}$ e a $MI_{P}$. 
In base ai valori forniti, si ha $MI_N <MI_E<MI_P$ (i pedici $E$, $P$, $N$ stanno evidentemente per \emph{equabile}, \emph{pitagorico}, \emph{naturale}). 

\noindent
L'intervallo $I$ di ampiezza $5/4$, prossimo a quello formato da $4$ semitoni equabili, verr\`a detto di \emph{terza maggiore naturale} (da cui il nome del Principio appena formulato) ed \`e caratterizzato da
\begin{equation}
\label{terzanat}
I(f_{i_1}, f_{i_2})=\frac{f_{i_1}}{f_{i_2}}=\frac{5}{4}
\end{equation}
\`E proprio la semplificazione del rapporto $81/64$ in favore di $5/4$ per l'intervallo di terza (di importanza fondamentale nella musica dal XVI secolo in poi) che ha comportato le novit\`a sulla scala pitagorica, riassumibili nell'importazione del numero $5$.

\noindent
La procedura illustrata per la scala pitagorica viene ripercorsa includendo dunque un ulteriore fattore primo rispetto alla (\ref{scalapit}): va ora considerato l'insieme dei suoni generati da $f_0$ 
\begin{equation}
\label{scalanat}
f=2^{m} 3^{n}5^{p}f_0, \qquad m, n, p \in {\Bbb Z}
\end{equation}
dove i numeri interi $m$, $n$ e $p$ possono essere positivi o negativi. 
Alle modalit\`a operative $(i)$--$(iii)$ della scala pitagorica, elencate appena dopo (\ref{scalapit}), vanno aggiunte le seguenti: 
\begin{itemize}
\item[$(iv)$] ogni moltiplicazione per $5$ porta il suono alla terza maggiore di due ottave sopra, ovvero lo trasferisce due ottave ed una terza maggiore sopra; l'esponente $p$ positivo fa eseguire $p$ di queste operazioni in avanti;
\item[$(v)$] ogni divisione per $5$ porta il suono due ottave ed una terza sotto; l'esponente $p$ negativo fa ripetere $-p$ volte l'operazione all'indietro, ovvero verso il grave. 
\end{itemize}
Le due procedure sono schematizzate come segue: a sinistra si passa da $f_0$ a $5f_0$, a destra da $f_0$ a $\frac{1}{5}f_0$. 

$$
\begin{array}{l}
\overbrace{|----|}^{ottava}--------\overbrace{|--|}^{terza\;maggiore} \\
f_0 \;\;\qquad \underbrace{2f_0\;\; \qquad \quad \;\;\;\qquad \quad \;\; 4f_0}_{ottava}\quad 5f_0
\end{array}\quad 
\begin{array}{l}
\overbrace{|--|}^{terza\;maggiore}------\overbrace{|------------|}^{ottava}\\
\quad \;\dfrac{1}{5}f_0\;\,
\underbrace{\dfrac{1}{4}f_0\; \;\,\;\;\quad \quad \quad \;\;\;\;\; \;\;\dfrac{1}{2}f_0}_{ottava}  \hspace{3.7truecm} f_0
\end{array}
$$

\noindent
Adoperando i medesimi argomenti utilizzati per la scala pitagorica (\ref{scalapit}), si dimostra che, per  infinit di $n, m, p\in {\Bbb Z}$, i numeri (\ref{scalanat}) sono razionali, sono infiniti, sono distinti se gli esponenti sono distinti.
In particolare, se $m\not =0$ e $p\not =0$ non si pu\`o ottenere una potenza del 2, ovvero a partire da $f_0$, nessun numero (\ref{scalanat}) si sovrapporr\`a pi\`u ad un'ottava di $f_0$.
Fissato un suono di riferimento $f_0$, chiamiamo \emph{scala naturale} o \emph{scala} \fbox{P} l'insieme dei suoni (\ref{scalanat}), infiniti e distinti. 

\subsection{La divisione armonica}

\noindent
Siamo nuovamente di fronte al problema di stabilire, a partire da $f_0$, una selezione in numero finito di suoni naturali (\ref{scalanat}), sempre al fine di rendere pratica la possibilit\`a di utilizzare la scala sugli strumenti musicali.

\noindent
Il procedimento svolto per la scala pitagorica, di trasferimenti in avanti e all'indietro a intervalli di quinte, fermandosi dopo un prefissato numero di passi, pu\`o essere riprodotto anche in questo caso, aggiungendo la possibilit\`a, come si \`e precisato, di trasferimenti in avanti e all'indietro a intervalli di terze.
Come si pu\`o immaginare, la situazione \`e ora molto pi\`u articolata e l'incasellatura dei suoni in quelli gi\`a scoperti sulle scale \fbox{E} e \fbox{P} potrebbe apparire in un certo senso artificiale.

\noindent
Se la procedura aritmetica della scala \fbox{N} si prefigura troppo articolata e contorta, 
vi \`e d'altra parte un metodo di tipo geometrico che compie la selezione dei suoni (\ref{scalanat}) in modo decisamente naturale: il metodo si basa su una costruzione geometrica spesso indicata come perfetto che connubio tra musica e architettura dal Rinascimento in poi e nota come 
\emph{divisione armonica}. La nozione di base \`e la seguente: 
sia $AB$ un segmento, $C$ un punto interno ad esso e $D$ un punto esterno sulla medesima retta che contiene il segmento:
$$
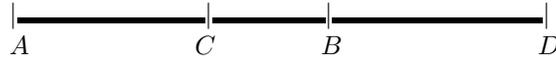
\begin{array}{l}
|\rule{25mm}{.7mm}|\rule{15mm}{.7mm}|\rule{28mm}{.7mm}|\\
A\hspace{2.2cm}C\hspace{1.4cm}B\hspace{2.6cm}D
\end{array}
$$
La quaterna $(A,B,C,D)$ forma una \emph{divisione armonica} se vale l'uguaglianza fra i rapporti
\begin{equation}
\label{proparm}
\dfrac{AC}{CB}=\dfrac{AD}{BD}
\end{equation}
Notiamo che all'avvicinarsi di $C$ verso $B$, anche il punto $D$ si avvicina a $B$, affinch\'e rimanga vera la (\ref{proparm}). D'altra parte, se $C$ si avvicina verso il punto medio del segmento $AB$ (dunque $AC/AB$ tende ad essere $1$), il punto $D$ deve allontanarsi all'infinito, per far s\`\i$\,$ che anche il rapporto $AD/BD$ sia sempre pi\`u prossimo a $1$.
Osserviamo inoltre che, essendo la definizione (\ref{proparm}) formulata sui rapporti, se $AC$ e $AD$ vengono moltiplicati per il medesimo numero $\gamma$, anche le suddivisioni $AB$ e $BD$ vengono scalate del medesimo fattore (ovvero: le distanze tra i punti $A$, $B$, $C$ $D$ si mantengono in scala).

\noindent
{\bf Proposizione}.  Se i punti $(A,B,C,D)$ formano una divisione armonica, allora 
\begin{equation}
\label{mediaarm}
\frac{1}{AB}=\frac{1}{2}\left(\frac{1}{AC}+\frac{1}{AD}\right)
\end{equation}

\noindent
{\bf Dim}. La (\ref{proparm}) equivale a $\frac{BC}{AC}=\frac{BD}{AD}$; sostituendo $BC=AB-AC$ e $BD=AD-AB$ si trova
$$
\frac{AB-AC}{AC}=\frac{AD-AB}{AD} \;\;\;\textrm{ovvero}\;\;\;\frac{AB}{AC}-1=1-\frac{AB}{AD}
$$
Dividendo per $AB$ e spostando i termini: $\frac{2}{AB}=\dfrac{1}{AC}+\frac{1}{AD}$ che equivale alla (\ref{mediaarm}).$\quad\square$

\noindent
La quantit\`a (\ref{mediaarm}) \`e collegata alla \emph{media armonica} di due numeri positivi $a$ e $b$, definita come il numero $m_H(a,b)=\frac{2}{\frac{1}{a}+\frac{1}{b}}$, 
ovvero la media armonica \`e il reciproco della media aritmetica fra i reciproci dei numeri dati $a$ e $b$.
Confrontando quest'ultima definizione con (\ref{mediaarm}), si  vede che
$AB=\frac{2}{\frac{1}{AC}+\frac{1}{AD}}=m_H(AC, AD)$, 
ossia $AB$ \`e la media armonica dei segmenti $AC$ e $AD$.

\noindent
Per i medesimi numeri $a$ e $b$ conosciamo la \emph{media aritmetica} $m_A(a,b)=\frac{a+b}{2}$ 
e la \emph{media geometrica} $m_G=\sqrt{ab}$:  
si pu\`o dimostrare che la media geometrica \`e medio proporzionale tra la media aritmetica e quella armonica:
$m_A:m_G=m_G:m_H$, 
pertanto valgono le disuguaglianze $m_H\leq m_G\leq m_A$, dato che $m_A\geq m_G$. 

\noindent
Le tre medie possono essere evidenziate graficamente con la seguente costruzione geometrica: si dispongono in modo adiacente due segmenti $AC$ e $CB$ di lunghezza rispettivamente $a$ e $b$ e si traccia la semicirconferenza di diametro $AB$. Il centro $O$ viene unito con $H$, punto sulla semicirconferenza di intersezione con la retta ortogonale ad $AB$ e passante per $C$. Infine, il punto $K$ \`e l'intersezione della retta per $C$ ortogonale al raggio 
$OH$ con il raggio medesimo. 

{\footnotesize 
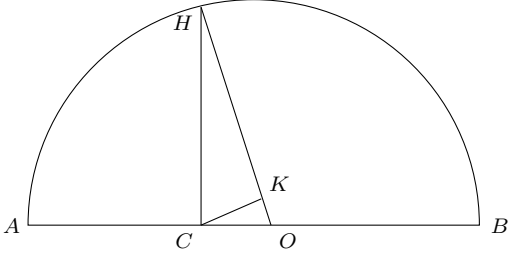
\begin{tikzpicture}
\draw[] (0,0) arc (180:0:3) -- cycle;

\coordinate [label=below left:\textcolor{black}{$A$}](a) at (0,0.2);
\coordinate [label=below left:\textcolor{black}{$B$}](b) at (6.5,0.2);
\coordinate [label=below left:\textcolor{black}{$C$}](c) at (2.3,0);
\coordinate [label=below right:\textcolor{black}{$O$}](o) at (3.23,0);
\coordinate [label=below left:\textcolor{black}{$H$}](h) at (2.3,2.9);
\coordinate [label=above right:\textcolor{black}{$K$}](k) at (3.10,0.35);
\draw [] (c) -- (h);
\draw [] (o) -- (h);
\draw [] (c) -- (k);
\end{tikzpicture}
}

\noindent 
I segmenti $OH$, $CH$ e $HK$ sono rispettivamente media aritmetica, media geometrica e media armonica di 
$AC=a$ e $CB=b$:
$$
m_A(a,b)=OH, \quad M_G(a,b)=CH, \quad m_H(a,b)=CK.
$$

\subsection{La selezione di suoni mediante la media armonica}

\noindent
L'idea che vogliamo sviluppare consiste nel compiere la scelta dei suoni (\ref{scalanat}) in modo che i gradini della scala si dispongano secondo divisioni armoniche: \`e come pensare ad una costruzione, i cui elementi divisori devono essere disposti in modo da dare origine a proporzioni armoniche nel numero maggiore possibile.
Pensiamo ed al segmento appena sopra la (\ref{proparm}) come una corda $AD$, in cui, a partire dalla suddivione effettuata in $C$, si voglia trovare il punto $B$ di divisione armonica: l'elemento favorevole consiste nel fatto che le frequenze sono inversamente proporzionali alle lunghezze delle corde, ovvero se $f_{AD}$ \`e la frequenza associata alla corda $AD$ in vibrazione, si ha
$f_{AD}=\frac{\kappa}{AD}$ con $\kappa$ costante e, analogamente, $f_{AB}=\frac{\kappa}{AB}$, $f_{AC}=\frac{\kappa}{AC}$, con la medesima costante $\kappa$. 
Sostituendo in (\ref{mediaarm}), vediamo che la divisione armonica $(A,B,C,D)$ comporta per le frequenze
\begin{equation}
\label{freqarm}
f_{AB}=\dfrac{1}{2}\left(f_{AC}+f_{AD}\right)
\end{equation}
ovvero la seguente

\noindent
{\bf Propriet\`a}. 
Se $B$ divide la corda $AD$ in modo che $AC$ e $AD$ formino la proporzione armonica (\ref{proparm}), allora la frequenza corrispondente alla corda $AB$ \`e la media aritmetica delle frequenze abbinate ad $AC$ e $AD$.

\noindent
Facendo appello al solo Principio dell'ottava e cercando il medio armonico $B$ di (\ref{proparm}), quando $AD=2AC$ (corda divisa a met\`a): da (\ref{proparm}), calcoliamo
$AB=\dfrac{2}{\dfrac{2}{AD}+\dfrac{1}{AD}}=\dfrac{2}{3}AD$
riscontrando che $B$ \`e a due terzi della lunghezza complessiva e produce il suono (\ref{freqarm}), se $f_0$ corrisponde alla corda intera $AD$: pertanto (\ref{freqarm}), che si scrive 
$f_{AB}=\dfrac{1}{2}(2f_0+f_0)=\dfrac{3}{2}f_0$
d\`a luogo all'intervallo di \emph{quinta pitagorica} (\ref{quintapit}). si \`e trovato che la corda intera e la corda divisa a met\`a hanno come medio armonico la corda che produce l'intervallo di quinta.
Per proseguire, cerchiamo la media armonica tra la corda intera e la corda di lunghezza $2/3$ rispetto ad essa: utilizzando ora la (\ref{proparm}) con $AC=\dfrac{2}{3}AD$, si trova (continuando a chiamare $AB$ il segmento trovato)
$AB=\dfrac{2}{\dfrac{3}{2AD}+\dfrac{1}{AD}}=\dfrac{4}{5}AD$
con frequenza (\ref{freqarm}) pari a $f_{AB}=\dfrac{5}{4}f_0$, corrispondente all'intervallo di \emph{terza naturale}) (\ref{terzanat}). Applicando ancora una volta la procedura, si calcola la media armonica con $AC=\dfrac{5}{4}AD$, per trovare
$AB=\dfrac{2}{\dfrac{5}{4AD}+\dfrac{1}{AD}}=\dfrac{8}{9}AD$
con la frequenza $f_0=\dfrac{9}{8}f_0$ che d\`a luogo all'intervallo coincidente con il suono pitagorico $\dfrac{3^2}{2^3}f_0$ che, per la vicinanza al tono equabile $\sqrt[12]{2^2}$, chiameremo \emph{tono pitagorico}. 

\noindent
Conviene riassumere queste significative scoperte di disposizioni armoniche e i corrispondenti esiti in termini di suoni nel seguenta schema, in cui la doppia stanghetta segna la divisione armonica in $B$ e l'intervallo sottolineato \`e quello prodotto per proporzione armonica dall'intervallo non sottolineato: 
{\footnotesize 
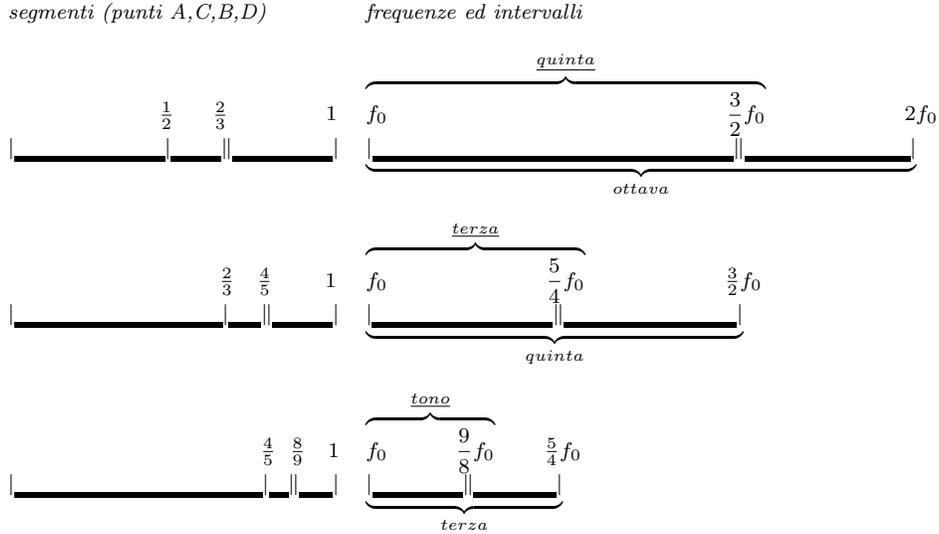
$$
\begin{array}{ll}
 \emph{segmenti (punti A,C,B,D)} & \emph{frequenze ed intervalli}\\
\\
\hspace{2cm}\frac{1}{2} \hspace{.5cm} \frac{2}{3} \hspace{1.3cm}1  &\overbrace{f_0 \hspace{4.5cm} 
\frac{3}{2}f_0}^{\underline{quinta}} \hspace{1.8cm} 2f_0 \\
|\rule[-1mm]{2cm}{.7mm}|\rule[-1mm]{.6666cm}{.7mm}\|\rule[-1mm]{1.3333cm}{.7mm}|& 
\underbrace{|\rule[-1mm]{4.8cm}{.7mm}\|\rule[-1mm]{2.2cm}{.7mm}|}_{ottava} \\
\\
 \hspace{2.8cm}\frac{2}{3} \hspace{.29cm} \frac{4}{5}\hspace{.7cm}1 & \overbrace{f_0 \hspace{2.1cm}\frac{5}{4}f_0}^{\underline{terza}} \hspace{1.8cm}\frac{3}{2}f_0  \\
|\rule[-1mm]{2.7666cm}{.7mm}|\rule[-1mm]{.4333cm}{.7mm}\|\rule[-1mm]{.8cm}{.7mm}| & 
\underbrace{|\rule[-1mm]{2.4cm}{.7mm}\|\rule[-1mm]{2.3cm}{.7mm}|}_{quinta} \\
\\
\hspace{3.35cm}\frac{4}{5} \hspace{.18cm} \frac{8}{9}\hspace{.3cm}1 & \overbrace{f_0\hspace{.9cm} 
\frac{9}{8}f_0}^{\underline{tono}} \hspace{.6cm}\frac{5}{4}f_0 \\
|\rule[-1mm]{3.3cm}{.7mm}|\rule[-1mm]{.2555cm}{.7mm}\|\rule[-1mm]{.4444cm}{.7mm}| & 
\underbrace{|\rule[-1mm]{1.2cm}{.7mm}\|\rule[-1mm]{1.1cm}{.7mm}|}_{terza}
\end{array}
$$
}
Impiegando le proporzioni armoniche, abbiamo ottenuto per il momento la scala di suoni armonici 
$$
f_0=DO, \quad \frac{9}{8}f_0=RE, \quad \frac{5}{4}f_0=MI, \quad \frac{3}{2} f_0=SOL, \quad 2f_0=DO
$$
Se si proseguisse con il calcolo delle medie armoniche utilizzando l'ultimo suono ottenuto $\frac{9}{8}f_0$, si otterrebbe la frequenza $\frac{17}{16}f_0$ che non fa parte di (\ref{scalanat}). Allo stesso modo, \`e facile constatare che le medie armoniche dei tre suoni trovati 
$\frac{3}{2}f_0$, $\frac{5}{4}f_0$ e $\frac{9}{8}f_0$ sono suoni estranei alla scala naturale (\ref{scalanat}): ad esempio, $\frac{1}{2}\left(\frac{5}{4}f_0+\frac{3}{2}f_0\right)= \frac{11}{8}f_0$.

\noindent
Possiamo tuttavia formulare un argomento che permette di completare la scala imponendo delle proporzioni armoniche tra i suoni mancanti e quelli gi\`a presenti: 
valutando gli spazi che lasciano i cinque suoni gi\`a trovati fra loro, diciamo di voler aggiungere, 
contando anche sull'esperienza gi\`a compiuta sulle scale precedenti, un suono tra $\frac{5}{4}f_0$ e $\frac{3}{2}f_0$ e due suoni tra $\frac{3}{2}f_0$ e $2f_0$. Iniziamo col dimostrare la seguente

\noindent
{\bf Proposizione}. 
Esistono in modo unico due suoni naturali $f_1$ e $f_2$  che formano proporzioni armoniche con il suono di base $f_0$ e l'ottava $2f_0$:
\begin{equation}
\label{quartasesta}
f_1=\dfrac{1}{2}\left(f_0+f_2\right), \quad f_2=\dfrac{1}{2}\left(f_1+2f_0\right)
\end{equation}

\noindent
{\bf Dim}. Sviluppando la (\ref{quartasesta}) troviamo
$2f_1=f_0+f_2$, $2f_1=4f_2-4f_0$
da cui, per confronto, si ottiene $f_2=\frac{5}{3}f_0$ e per sostituzione $f_1=\frac{4}{3}f_0$. 
$\quad\square$

\noindent
I due suoni ottenuti, l'uno tra la terza e la quinta $\frac{5}{4}f_0<f_1<\frac{3}{2}f_0$, l'altro tra la quinta e l'ottava $\frac{3}{2}f_0<f_2<2f_0$, vanno a definire rispettivamente il \emph{FA naturale}
$FA_N=\frac{4}{3}f_0$ ed il \emph{LA naturale} $LA_N=\frac{5}{3}f_0$. 
I corrispondenti intervalli $I\left(f_0, \frac{4}{3}f_0\right)=\frac{4}{3}=1,{\bar 3}$ e $I\left(f_0, \frac{5}{3}f_0\right)=\frac{5}{3}=1,{\bar 6}$ sono gli intervalli di \emph{quarta} e di \emph{sesta naturali}. 
Rispetto alle scale precedenti, valgolo le relazioni
$$
FA_N=\dfrac{4}{3}f_0=FA_P<FA_E=\sqrt[12]{2^5}f_0, \qquad
LA_N=\dfrac{5}{3}f_0<LA_E=\sqrt[12]{2^9}f_0  < LA_P=\dfrac{27}{16}f_0
$$
dove per il suono $LA$ la differenza non \`e certamente insignificante. Lo ``spazio'' lasciato fra $\frac{5}{3}f_0$ e $2f_0$ impone il completamento della scala naturale: dimostriamo un risultato che produrr\`a univocamente la selezione del suono mancante, sulla base della disposizione armonica dei suoni gi\`a presenti.

\noindent
{\bf Proposizione}. Esiste un'unica possibilit\`a di formare una proporzione armonica tra due suoni naturali scelti nell'insieme ${\cal I}=(f_0, \frac{9}{8}f_0, \frac{4}{3}f_0, \frac{5}{4}f_0, \frac{3}{2}f_0, \frac{5}{3}f_0, 2f_0)$ con un suono naturale $f_3$ compreso tra $\frac{5}{3}f_0$ e $2f_0$.

\noindent
{\bf Dim}.  La richiesta \`e
$f_{N_1}=\frac{1}{2}\left(f_{N_2}+f_3\right)$ con  
$f_{N_1}$, $f_{N_2}\in {\cal I}$ e la restrizione $\frac{5}{3}f_0<f_3<2f_0$. Ricaviamo
$f_3=2f_{N_1}-f_{N_2}$ a cui va imposto $\frac{5}{3}f_0<2f_{N_1}-f_{N_2}<2f_0$:
esaminando tutte le possibilit\`a, si riscontra facilmente che $N_1$ pu\`o assumere solo i valori $\frac{5}{3}f_0$ e $\frac{3}{2}f_0$. 
Per $f_{N_1}=\frac{5}{3}f_0$ gli unici due suoni $f_{N_2}$ compatibili con la restrizione sono $f_{N_2}=\frac{5}{4}f_0, \frac{3}{2}f_0$, i quali per\`o darebbero luogo a $f_3=\frac{35}{12}f_0, \frac{11}{6}f_0$ che non fanno parte dei numeri (\ref{scalanat}).   
Per quanto riguarda invece $f_{N_1}=\frac{3}{2}f_0$, le uniche possibilit\`a per $N_2$ corrispondono a  $f_{N_2}=\frac{9}{8}f_0$ e $\frac{5}{4}f_0$, che danno origine rispettivamente a $f_3=\frac{15}{8}f_0$ e $\frac{7}{4}f_0$, quest'ultimo suono non appartenente alla scala (\ref{scalanat}): si conclude dunque che l'unica circostanza compatibilie con le richieste \`e
$f_{N_1}=\frac{3}{2}f_0$, $f_{N_2}=\frac{9}{8}f_0$, che conducono a $f_3=2\frac{3}{2}f_0-\frac{9}{8}f_0=\frac{15}{8}f_0$ $\quad\square$. 

\noindent
Il suono trovato \`e il \emph{SI della scala naturale} $SI_N$, compreso tra LA=$\frac{5}{3}f_0$ e DO=$2f_0$. 
Il confronto con i suoni omologhi delle altre scale porta a 
$$
SI_N=\frac{15}{8}f_0<SI_E=\sqrt[12]{2^{11}}f_0<SI_P=\frac{243}{128}f_0.
$$
L'intervallo $I\left(f_0, \dfrac{15}{8}f_0\right)=\dfrac{15}{8}$ \`e l'\emph{intervallo di settima} della scala naturale.

\subsection{La scala naturale diatonica}

\noindent
Il procedimento esposto ha individuato 6 suoni intermedi del tipo (\ref{scalanat}). Riassumiamo la procedura, sottolineando l'utilizzo della proporzione armonica come mezzo eslusivo per generare i suoni 
e l'unicit\`a, in un certo senso, della scala ottenuta. Indichiamo con ${\cal M}_A(suono_1, suono_2)=suono_3$ la circostanza in cui $suono_1$ e $suono_2$ producono per proporzione armonica $suono_3$, ovvero $f_3=\frac{1}{2}\left(f_1+f_2\right)$ (come in (\ref{freqarm})), se $f_1$, $f_2$, $f_3$ sono le frequenze corrispondenti ai tre suoni. Nella sintesi che segue i nomi delle note si riferiscono alla scala naturale.
\begin{itemize}
\item[$\bullet$]
A partire da $f_0=DO$ e $2f_0=DO$ all'ottava sopra si costruiscono
${\cal M}_A(DO,2DO)=SOL=\frac{3}{2}f_0$, ${\cal M}_A\left(DO, SOL\right)=MI=\frac{5}{4}f_0$ e ${\cal M}_A(DO, MI)=RE=\dfrac{9}{8}f_0$;
\item[$\bullet$] dai tre suoni ottenuti non si ottengono per proporzione armonica altri suoni di tipo (\ref{scalanat});
\item[$\bullet$] esistono unicamente due suoni $FA$, $LA$ per cui 
$FA={\cal M}_A(DO, LA), \quad LA={\cal M}_A(FA, 2DO)$
e risultano di tipo (\ref{scalanat}): $FA=\dfrac{4}{3}f_0$, $LA=\dfrac{5}{3}f_0$; 
\item[$\bullet$] l'unico suono di tipo (\ref{scalanat}) compreso fra $LA$ e $2DO$ e prodotto da una media armonica di suoni fin qui trovati \`e $SI$ che verifica
$SOL={\cal M}_A(RE, SI)$.
\end{itemize}

\noindent
Nella grafica che segue, si vede che ciascun suono della scala naturale \`e interessato almeno una volta nei rapporti armonici con gli altri suoni; oltre alle tre proporzioni armoniche delineate nel Paragrafo 5.2 aggiungiamo 
{\footnotesize 
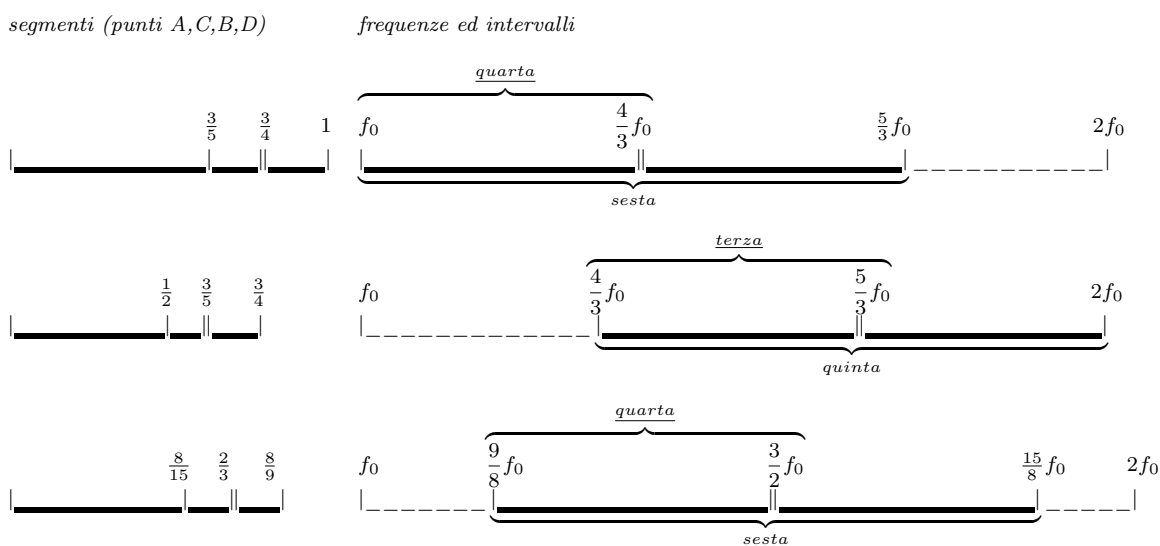
$$
\begin{array}{ll}
 \emph{segmenti (punti A,C,B,D)} & \emph{frequenze ed intervalli}\\
\\
\hspace{2.6cm}\frac{3}{5} \hspace{.5cm} \frac{3}{4} \hspace{.63cm}1 &
\overbrace{f_0 \hspace{3.1cm} \frac{4}{3}f_0}^{\underline{quarta}} \hspace{2.9cm} \frac{5}{3}f_0 \hspace{2.4cm} 2f_0\\
|\rule[-1mm]{2.55cm}{.7mm}|\rule[-1mm]{.6cm}{.7mm}\|\rule[-1mm]{.75cm}{.7mm}|& 
\underbrace{|\rule[-1mm]{3.6cm}{.7mm}\|\rule[-1mm]{3.4cm}{.7mm}|}_{sesta} \_\_\_\_\_\_\_\_\_\_\_|\\
\\
\hspace{2cm}\frac{1}{2} \hspace{.32cm}\frac{3}{5} \hspace{.48cm} \frac{3}{4} &
f_0 \hspace{2.7cm} \overbrace{ \frac{4}{3}f_0 \hspace{3cm}\frac{5}{3}f_0}^{\underline{terza}} \hspace{2.6cm}2f_0  \\
|\rule[-1mm]{2cm}{.7mm}|\rule[-1mm]{.4cm}{.7mm}\|\rule[-1mm]{.6cm}{.7mm}| & 
|\_\_\_\_\_\_\_\_\_\_\_\_\_ \underbrace{|\rule[-1mm]{3.35cm}{.7mm}\|\rule[-1mm]{3.15cm}{.7mm}|}_{quinta}\\
\\
\hspace{2.1cm} \frac{8}{15} \hspace{.32cm}\frac{2}{3} \hspace{.4cm} \frac{8}{9} & 
f_0\hspace{1.35cm} \overbrace{\frac{9}{8}f_0\hspace{3.2cm} \frac{3}{2}f_0}^{\underline{quarta}} \hspace{2.8cm}\frac{15}{8}f_0 \hspace{.8cm}2f_0 \\
|\rule[-1mm]{2.2333cm}{.7mm}|\rule[-1mm]{.5333cm}{.7mm}\|\rule[-1mm]{.5333cm}{.7mm}| & 
|\_\_\_\_\_\_\_\underbrace{|\rule[-1mm]{3.6cm}{.7mm}\|\rule[-1mm]{3.4cm}{.7mm}|}_{sesta}\_\_\_\_\_|
\end{array}
$$
}

\noindent
Il seguente schema riporta i suoni della scala naturale diatonica e gli intervalli che intercorrono fra essi:

$$
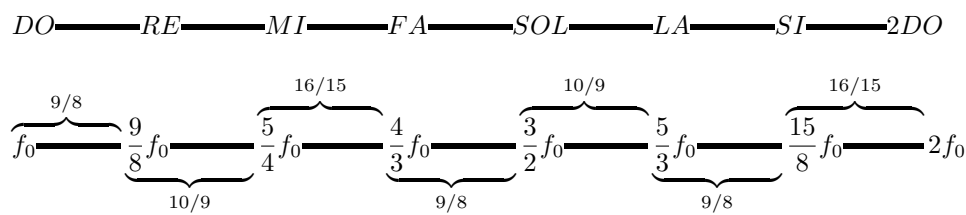
\begin{array}{l}
DO\rule[.5mm]{1.125cm}{.7mm}
RE \rule[.5mm]{1.111cm}{.7mm}
MI\rule[.5mm]{1.066cm}{.7mm}
FA\rule[.5mm]{1.125cm}{.7mm} 
SOL\rule[.5mm]{1.111cm}{.7mm} 
LA\rule[.5mm]{1.125cm}{.7mm}
SI\rule[.5mm]{1.066cm}{.7mm} 2DO\\
\\
\overbrace{f_0 \rule[.5mm]{1.125cm}{.7mm}}^{9/8}
\underbrace{\dfrac{9}{8}f_0 \rule[.5mm]{1.111cm}{.7mm}}_{10/9}
\overbrace{\dfrac{5}{4}f_0\rule[.5mm]{1.066cm}{.7mm}}^{16/15} 
\underbrace{\dfrac{4}{3}f_0\rule[.5mm]{1.125cm}{.7mm}}_{9/8} 
\overbrace{\dfrac{3}{2}f_0\rule[.5mm]{1.111cm}{.7mm}}^{10/9} 
\underbrace{\dfrac{5}{3}f_0\rule[.5mm]{1.125cm}{.7mm}}_{9/8} 
\overbrace{\dfrac{15}{8}f_0\rule[.5mm]{1.066cm}{.7mm}}^{16/15} 2f_0
\end{array}
$$

\noindent
La scala naturale (non solo diatonica, ma comprendente i suoni cromatici) fu codificata in importanti trattati musicali come \cite{zar} ci si riferisce ad essa anche con il termine di \emph{scala zarliniana}. 
L'evidente vantaggio consiste nella maggiore semplicit\`a dei rapporti fra i numeri razionali che formano i suoni: ad esempio, $SI_P=(243/128)f_0$, mentre $SI_N=(15/8)f_0$. 
L'inconveniente sta per\`o nel fatto che i passaggi da un tono all'altro avvengono con due differenti rapporti, ovvero $9/8$ e $10/9$: questo rende la costruzione della scala fortemente dipendente dal suono iniziale $f_0$ ed impossibile la propriet\`a di trasposizione. Dal punto di vista musicale, esisteranno tonalit\`a ``privilegiate'' da scale pi\`u semplici, a scapito di altre, con rapporti pi\`u complessi. Al tempo stesso, come gi\`a avviene sulla scala pitagorica, alcuni intervalli visti identici dalla scala equabile, sono di ampiezza differente: ad esempio $I(DO,SOL)=3/2=1,5$, $I(RE,LA)=40/27=1,{\overline 148}$. L'inconveniente \`e particolarmente rilevante sugli strumenti ad intonazione fissa, come quelli a tastiera.
In effetti, con il sopravvento dal XVI secolo in poi della musica strumentale, si apre il vasto dibattito fra i toni della musica teorica e quelli impiegati dalla prassi musicale. La soluzione proposta dalla scala equabile \`e da collocarsi in realt\`a in una fase decisamente matura del dibattito: per la sua spontaneit\`a matematica e per le sue propriet\`a di omogeneit\`a \`e stato qui deciso di presentarla come prima.

\noindent
In conclusione replichiamo, limitatamente ai suoni diatonici, i valori normalizzati via via ottenuti per le tre scale; da essi si pu\`o dedurre quali sono i suoni ``calanti'' o ``crescenti'' rispetto agli omologhi e la lunghezza dei vari intervalli, questi ultimi pure confrontabili (ad esempio, si vede che la quinta equabile \`e pi\`u corta rispetto a quella pitagorica, coincidente con quella naturale, oppure che l'intervallo di terza \`e sensibilmente diverso nelle tre scale). Si mette in evidenza anche la struttura di tipo (\ref{ognifreq}), (\ref{scalapit}) e (\ref{scalanat}), ripettivamente per la scala $\fbox{E}$, $\fbox{P}$ e $\fbox{N}$.

\begin{center}
\begin{tabular}{|l||l|l | l ||}\hline
& \fbox{E} & \fbox{P} & \fbox{N} \\ \hline 
DO & $1$ & $1$ & $1$  \\
RE    & $\sqrt[12]{2^2}=1,12246...$ & $2^{-3}3^2=1,125$ & $2^{-3}3^2=1,125$ \\
MI   & $\sqrt[12]{2^4}=1,25992...$ & $2^{-6}3^4=1,26562...$ &  $2^{-2}5=1,25$ \\
FA    & $\sqrt[12]{2^5}=1,33484...$ & $2^23^{-1}=1,{\bar 3}$ & $2^23^{-1}=1,{\bar 3}$ \\
SOL   & $\sqrt[12]{2^7}=1,49830...$ & $2^{-1}3=1,5$ & $2^{-1}3=1,5$ \\
LA    & $\sqrt[12]{2^9}=1,68180...$ & $2^{-4}3^3=1,6875$& $3^{-1}5=1,{\bar 6}$\\
SI      & $\sqrt[12]{2^{11}}=1,88775...$ & $ 2^{-7}3^5 = 1,89843...$& $2^{-3}3 \times 5=1,875$ \\   
DO     & $2$  & $2$  & $2$ \\   
\hline
\end{tabular}
\end{center}

\section{Coda} 

\noindent
Concedendo un privilegio ad una prospettiva, per cos\`\i$\,$dire, assiomatica, si pu\`o operare questa sintesi: la scala equabile necessita solo del principio dell'ottava (collegabile al numero due) e di un concetto di distanza, quella pitagorica estende il principio accettando il numero tre, infine la scala naturale lo amplia ulteriormente incorporando il cinque.

\noindent
In modo forse spietato ed irriverente verso le tante tabelle sopra esibite di suoni in antagonismo tra loro per qualche decimale, viene in mente il celebre incontro tra J.~Brahms, non certo l'ultimo tra i musicisti\footnote{A detta dell'autore di questo articolo, il primo.} e lo straodinario scienziato H.~v.~Helmholtz\footnote{A detta di chiunque il fondatore dell'acustica come scienza moderna, dalla natura fisica del suono all'aspetto della sensazione percettiva dello stesso, vedi \cite{hel}. Fra i tanti progressi del fisico tedesco va ricordato il modo di rendere quantitativo il colore del suono (timbro) attraverso la scomposizione in suoni armonici; l'opera di Helmholtz \`e la ripercussione 
in musica dello straordinario progresso scientifico del XVIII e XIX secolo, ad opera di D'Alembert, Lagrange, Fourier, J.~W.~S.~Rayleigh.}: quest'ultimo fu un anacronistico sostenitore dell'intonazione naturale ({\it Reine Stimmung}, corrispondente alla nostra scala $\fbox{N}$, a scapito dell'ormai dilagante temperamento equabile\footnote{ ...{\it la meccanica degli strumenti e la sua comoda praticabilit\`a minacciano di aver ragione delle esigenze naturali dell'orecchio}..., da \cite{hel}}$\fbox{E}$ e chiam\`o a testimonianza dell'accordatura pura il celebre violinista J\'ozsef Joachim\footnote{egli collabor\`o con Brahms nella composizione dell'unico concerto per violino scritto nel 1878, fornendo preziosi consigli dal punto di vista della tecnica violinistica.}, il quale, secondo Helmholtz, eseguiva ``{\it le terze secondo il rapporto} $4/5$''.
Il biografo Kalbeck, autore di una importante biografia di Bramhs, riporta di quest'ultimo la seguente cronaca (la seguente traduzione \`e tratta da \cite{ser}), in prima persona:

\begin{minipage}[c]{.93\textwidth}

{\it Joachim ed io eravamo da Helmholtz, che ci illustrava le sue scoperte e le armonie pure con uno strumento da lui ideato\footnote{l'harmonium di Bosanquet, che suddivide l'ottava in $53$ parti uguali, presentato all'Esposizione Universale di Londra del 1876}. affermava che la terza dovesse essere un p\`o pi\`u alta, la settima pi\`u bassa dell'usuale. Joachim, che \`e persona di squisita cortesia, asseriva di aver ricevuto dagli intervalli un'impressione del tutto particolare, e fingeva di udirli allo stesso modo di Helmholtz. Allora gli dissi che la cosa era troppo seria perch\'e le buone maniere dovessero avere il sopravvento [...]. In campo musicale Helmholtz \`e proprio un atroce dilettante.}

\end{minipage}

\noindent
L'episodio, senz'altro ameno ma al tempo stesso educativo nei confronti di chi troppo insiste sull'irrinunciabile ruolo delle frazioni di tono in musica, delle quali la storia ci ha tramandato una sterminata botanica di microintervalli e termini tecnici, non pu\`o comunque distogliere dalla consapevolezza che i numeri hanno da sempre avuto un ruolo speciale nell'arte dei suoni, sin dai primi tempi del Quadrivium delle arti liberali Aritmetica, Geometria, Astronomia e Musica, quest'ultima \`e definita come `disciplina che tratta i numeri.

\noindent
Il rapporto speciale con la matematica, rispetto  alle nozioni matematiche diffusamente applicate alla comprensione dei fenomeni naturali e delle attivit\`a umane, comprese le arti,  pu\`o forse essere giustificato dal fatto che la sfuggente struttura della musica \`e ricostruita nell'esecuzione e nell'ascolto, deve catturare qualcosa in movimento.

\noindent
La musica \`e la capacit\`a di recepire, discernere, organizzare e apprezzare l'ordine nella totalit\`a dei movimenti intorno a noi. \`E l'incontro delle due razionalit\`a, l'una predisposta dall'Universo l'altra costruita dalle regole del pensiero. Pi\`u che altrove, una sola mente non basta a cogliere il tutto, occorre la cultura musicale per tramandare tutto ci\`o che si \`e raccolto.

\noindent
In modo pi\`u evidente e stringente che in altre discipline artistiche, le nozioni della pratica si consolidano, si cristallizzano in una teoria logica. Ma \`e la musica che stabilisce le regole: il ribaltamento dei ruoli di un'imposizione di regole pianificate astrattamente e programmate a tavolino 
non porta lontano, come \`e successo per la Dodecafonia.

\noindent
L'avanzamento formale e l'evoluzione emotiva della musica necessita che questa vada scritta, vada codificata: intervengono elementi razionali, astratti, annoverati nell'universo della 
matematica. Gli stessi studi ed esperimenti in campo acustico, sin dall'antichit\`a, hanno elargito entusiasmo riguardo la presenza della matematica nella musica.

\noindent
Per questa esigenza c'\`e chi ha visto, ad esempio, nel tetragramma o nel pentagramma (il riferimento che in verticale d\`a la frequenza, in orizzontale il tempo, con unit\`a di misura verticale la chiave, con unit\`a di misura orizzontale la battuta) un'anticipazione del sistema di riferimento cartesiano.

\noindent
\`E proprio attorno all'evoluzione della notazione musicale che effettuiamo 
la nostra riflessione finale: il pessimismo del primo grande teorico musicale del Medioevo, Isidoro da Siviglia, che affermava che la musica non potr\`a mai essere scritta, \`e stato capovolto dalla prodigiosa evoluzione formale culminante nelle formidabili conquiste in epoca beethoveniana, di compiuta possibilit\`a di tradurre sulla carta ci\`o che il compositore ha in mente.
Il successo \`e indubbiamente di carattere tecnico, operativo, possiamo dire scientifico.

\noindent
Sulla traccia di tale riuscita, altri aspetti del binomio scienza--musica potrebbero essere afferrati: pu\`o essere, ci domandiamo, una questione di capacit\`a scientifica di esplorare le discipline artistiche a porre il limite della comprensione e a causare schieramenti spesso fortemente contrapposti? 
Sar\`a mai ribaltato il pessimismo\footnote{``{\it La matematica e la fisica non hanno sensibilit\`a 
musicale}'', \cite{rig} p.~88} di un grande studioso di acustica o l'episodico scetticismo\footnote{Per il momento, lo spettro armonico prodotto da un oggetto lasciato cadere su un tasto del pianoforte e quello prodotto dal tocco di un famoso pianista sul medesimo tasto sono esattamente gli stessi.} 
in \cite{jea} di un illustre scienziato britannico dall'ampliamento della conoscenza scientifica?

\noindent
Non sono certo assenti argomenti e sforzi in questa direzione, ad esempio elaborati con l'Algebra, oppure con il formalismo esasperato delle strutture matematiche di \cite{xen}), oppure con la recente teoria degli insiemi fuzzy, per convogliare la nozione di intervallo in un insieme tollerato di frequenze (\cite{fuz1}), comprensivo delle sfumature con cui viene quantificato il medesimo intervallo in pi\`u scale, oppure per affrontare questioni legate alla percezione emotiva (\cite{fuz}).

\noindent
Potr\`a mai la scienza essere in grado di spiegare perch\'e, in musica, per imprimere la sensazione di andare a tempo non si deve andare dietro la scansione esatta del metronomo, oppure perch\'e le tonalit\`a, tutte uguali e indistinguibili secondo la scala $\fbox{E}$, hanno ciascuna un colore ed un carattere diverso dalle altre?

\noindent
L'autore di questo articolo immagina un p\`o in questo modo l'adoperarsi per procurare un ponteggio ma\-te\-ma\-ti\-co al mondo dei suoni, con il risultato di un capriccio, come genere dell'arte figurativa\footnote{Esempi mirabili sono i Capricci di G.~B.~Piranesi (1720--1778) e B.~Bellotto (1721--1780).}: non \`e la mancanza di forma a caratterizzare questi generi pi\`u liberi, bens\`i il fatto che una forma viene inventata, trovata strada facendo: un meraviglioso, sistematico ed elaborato accostamento di architetture, che tuttavia non esiste.


\end{document}